\documentclass[pdflatex,sn-mathphys-num]{sn-jnl}

\usepackage{graphicx}
\usepackage{multirow}
\usepackage{amsmath,amssymb,amsfonts}
\usepackage{amsthm}
\usepackage{mathrsfs}
\usepackage[title]{appendix}
\usepackage{xcolor}
\usepackage{textcomp}
\usepackage{manyfoot}
\usepackage{booktabs}
\usepackage{algorithm}
\usepackage{algorithmicx}
\usepackage{algpseudocode}
\usepackage{listings}

\usepackage{mathtools}
\usepackage{enumitem}
\usepackage{comment}

\numberwithin{equation}{section}

\theoremstyle{thmstylethree}
\newtheorem{definition}{Definition}[section]

\theoremstyle{thmstyleone}

\newtheorem{proposition}[definition]{Proposition}
\newtheorem{lemma}[definition]{Lemma}
\newtheorem{theorem}[definition]{Theorem}
\newtheorem{corollary}[definition]{Corollary}

\theoremstyle{thmstyletwo}
\newtheorem{remark}[definition]{Remark}

\newcommand{\R}{\mathbb R}
\newcommand{\C}{\mathbb C}

\newcommand{\Cp}{\mathbb C_+}
\newcommand{\Cm}{\mathbb C_-}
\newcommand{\E}{\mathbb E}

\newcommand{\dd}{\,d}
\newcommand{\inner}[2]{\left\langle #1,#2\right\rangle}
\newcommand{\cB}{\mathcal B_\alpha}

\newcommand{\1}{\mathbf 1}

\begin{document}

\title[Article Title]{Sparse Analytic Representation of Stationary Processes with Bounded Spectral Density}


\author[1]{\fnm{Yiqiao} \sur{Zhang}}
\nomail

\author[2]{\fnm{Fang} \sur{Yang}}
\nomail

\author*[3]{\fnm{Tao} \sur{Qian}}
\email{tqian@must.edu.mo}


\affil[1]{
	School of Computer Science and Engineering,
	Faculty of Innovation Engineering,
	Macau University of Science and Technology,
	Macao 999078, China
}

\affil[2]{
	Department of Physics,
	Zhejiang Normal University,
	Jinhua 321004, China
}

\affil[3]{
	Macau Center for Mathematical Sciences, School of Computer Science and Engineering, Macau University of Science and Technology, Taipa, 999078, Macao, China
}

\abstract { \unboldmath We consider centered second-order stationary processes	\(Y=\{Y_t\}_{t\in\mathbb R}\) with bounded spectral density \(S_Y\)
	and develop a sparse analytic representation framework. A nontrivial stationary process on the whole real line generally does not belong to \(L^2(\Omega,L^2(\mathbb R))\), since stationarity implies \(\mathbb E\|Y\|_{L^2(\mathbb R)}^2=\infty\). To overcome this infinite-energy obstruction, we replace the sample-path
	viewpoint by the bounded linear operator
	$$
	\widetilde Y:L^2(\mathbb R)\longrightarrow L^2(\Omega),
	\qquad
	\widetilde Y(f)
	=
	\int_{\mathbb R}\widehat f(\xi)\,Z_Y(\dd\xi),
	$$
	where \(\xi\) denotes the frequency variable and \(Z_Y\) is the orthogonal
	random measure appearing in the spectral representation of \(Y\).
	Applying \(\widetilde Y\) to Cauchy kernels in the upper and lower half-planes yields random analytic functions \(U_+\) and \(U_-\). 
	For \(\alpha>1\), we prove that
	$$
	U_+\in L^2\bigl(\Omega,B_\alpha^2(\mathbb C^+)\bigr),
	\qquad
	U_-\in L^2\bigl(\Omega,B_\alpha^2(\mathbb C^-)\bigr),
	$$
	where \(B_\alpha^2(\mathbb C^\pm)\) are weighted Bergman spaces on the upper and lower half-planes, respectively. We then combine \(U_+\) and \(U_-\) on \(\mathbb C^+\cup\mathbb C^-\) into the random analytic representation
	$$
	U(\omega,z)
	=
	\begin{cases}
		U_+(\omega,z), & z\in\mathbb C^+,\\[2pt]
		U_-(\omega,z), & z\in\mathbb C^-,
	\end{cases}
	$$
	and, with
	$
	\mathcal B_\alpha
	=
	B_\alpha^2(\mathbb C^+)\oplus B_\alpha^2(\mathbb C^-),
	$
	prove that
	$
	U\in L^2(\Omega,\mathcal B_\alpha).
	$
	We establish the
	boundary vanishing condition, energy identity, and convergence, and derive
	covariance identities for the random coefficients. 
	In the Bochner space $L^2(\Omega,\mathcal B_\alpha)$,
	we apply stochastic pre-orthogonal adaptive Fourier
	decomposition (SPOAFD) separately to $U_+$ and $U_-$,
	using the normalized reproducing kernels of the
	corresponding weighted Bergman spaces as dictionaries.
	The adaptive selection is determined by the spectral density of the original stationary process. Finally, using the Poisson
	kernel and the SPOAFD expansions on both half-planes, we establish
	the weak mean-square convergence of the original stationary process
	against test functions in \(C_c^\infty(\mathbb R)\).}


\maketitle

\section{Introduction}

In this section, we present the problem setting and outline our approach. More detailed analysis and proofs are provided in the subsequent sections.

Sparse representation methods have been widely developed for deterministic signals and operators, typically assumed to have finite energy in a Hilbert space such as $L^2$ or, for holomorphic functions, the Hardy space $H^2$. Among these methods, Adaptive Fourier Decomposition (AFD) provides adaptive rational approximations using parameterized Szeg? kernels and Takenaka--Malmquist systems. At each step, an energy matching principle selects the parameter, yielding an orthonormal system adapted to the signal \cite{QianWang2011,Qian2016}.

The deterministic AFD theory has been developed in several directions. Variants include weak AFD \cite{QianWang2013,Qian2015Book}, \(n\)-best AFD, associated with  \(n\)-best rational approximation \cite{Qian2014Cyclic}, and unwinding Blaschke expansions based on successive Blaschke factorizations \cite{CoifmanSteinerberger2017,CoifmanSteinerbergerWu2017}. Pre-Orthogonal Adaptive Fourier Decomposition (POAFD) extends the underlying adaptive selection principle from Hardy spaces to general Hilbert spaces and reproducing kernel Hilbert spaces equipped with dictionaries \cite{Qian2020}. More recently, Double-Zero AFD has been introduced by exploiting the additional interpolation property generated by the optimal parameter selection \cite{QianWuQuWang2026}. Related developments also include matrix-valued and multi-channel signals \cite{WangWongRosaQianWan2022,AlpayColomboQianSabadini2017}, as well as quaternionic and other hypercomplex settings \cite{QianSproessigWang2012}. These developments have made AFD-type methods applicable to a broad range of approximation and signal-analysis problems; an overview of this line of research can be found in \cite{Qian2018Survey}.

The deterministic theory was later extended to random signals. In
\cite{Qian2022}, AFD was extended to  the Hardy space-valued random signals,
leading to stochastic AFD (SAFD), and subsequently to stochastic
pre-orthogonal AFD (SPOAFD) in general Hilbert spaces with suitable
dictionaries. In RKHSs, normalized reproducing kernels provide a natural choice of dictionary for SPOAFD. SAFD and SPOAFD were further compared with the
Karhunen--Lo\`eve (KL) expansion in \cite{QianZhangLiuQu2023}. Under the
assumptions considered there, the AFD-type methods achieve an
\(O(n^{-1/2})\) approximation estimate comparable to that considered for KL,
while avoiding covariance eigendecomposition and allowing problem-oriented dictionaries, such as potential kernels.


The SAFD and SPOAFD theory described above is developed
for random signals in the Bochner space
$L^2(\Omega;\mathcal H)$
\cite{Qian2022,QianZhangLiuQu2023}.
This space is defined by (see \cite{vanNeerven2008})
\[
L^2(\Omega;\mathcal H)
:=
\left\{
F:\Omega\to\mathcal H:
F \text{ is an $\mathcal H$-valued random variable and }
\mathbb E\|F\|_{\mathcal H}^2<\infty
\right\},
\]
with norm
\(\|F\|_{L^2(\Omega;\mathcal H)}^2
=
\mathbb E\|F\|_{\mathcal H}^2.\)
This provides the Hilbert-space setting for measuring stochastic approximation errors.
We write
$$
L^2(\Omega):=L^2(\Omega,\mathbb C)
=
\left\{
X:\Omega\to\mathbb C:
X \text{ is a complex-valued random variable and }
\mathbb E|X|^2<\infty
\right\}.
$$
for the usual space of random variables with finite second moments.

A difficulty occurs when one tries to apply this framework directly to a stationary process defined on the entire real line. Let \(Y=\{Y_t:t\in\mathbb R\}\) be a non-trivial centered second-order stationary process. If
\[
\mathbb E Y_t=0,
\qquad
\mathbb E|Y_t|^2=\sigma^2>0,
\qquad
t\in\mathbb R,
\]
then, by stationarity and Tonelli's theorem,
\[
\mathbb E\int_{\mathbb R}|Y_t|^2\dd t
=
\int_{\mathbb R}\mathbb E|Y_t|^2\dd t
=
\int_{\mathbb R}\sigma^2\dd t
=
\infty.
\]
Consequently,
\(Y\notin L^2\bigl(\Omega,L^2(\mathbb R)\bigr).\)

Thus, a stationary process on the whole real line cannot be regarded as a random element of the Bochner space
\(L^2\bigl(\Omega,L^2(\mathbb R)\bigr)\).
This makes it difficult to apply the AFD methodology to stationary random processes in the usual Bochner-space setting.

For a general wide-sense stationary process, stationarity alone
does not guarantee that its sample paths are Lebesgue measurable
or belong to any particular classical function space.
In particular, membership in $L^1_{\mathrm{loc}}(\mathbb R)$
or $L^2_{\mathrm{loc}}(\mathbb R)$ cannot be concluded without
additional assumptions. If the process is jointly measurable,
then, by its finite second moments and wide-sense stationarity,
its sample paths belong to $L^2_{\mathrm{loc}}(\mathbb R)$,
and hence to $L^1_{\mathrm{loc}}(\mathbb R)$, almost surely.
Moreover, almost every sample path can be identified with a
tempered distribution in $\mathcal S'(\mathbb R)$.
These conclusions do not guarantee membership in the global
spaces $L^1(\mathbb R)$ or $L^2(\mathbb R)$, nor do they imply
pathwise continuity or differentiability.

The present study shows that, for a stationary process with
bounded spectral density, the associated generalized random
operator can be applied to the Cauchy kernels of the upper
and lower complex half-planes. The resulting random functions
admit analytic versions whose sample paths belong almost surely
to certain weighted Bergman spaces. We then construct SPOAFD sparse representations
of these random analytic functions, establish their
convergence, and recover the original stationary
process in the weak mean-square sense.

To establish these results, we describe the stationary process
through a bounded linear operator acting on deterministic
test functions.
Let \(Y=\{Y_t:t\in\mathbb R\}\) be a centered second-order stationary process with spectral density \(S_Y\) and spectral measure \(\mu_Y\) given by \eqref{eq:spectral-measure}. Then \(S_Y\in L^1(\mathbb R)\) and \(S_Y\ge 0\) almost everywhere. We further assume that \(S_Y\in L^\infty(\mathbb R)\).
By the spectral representation theorem, \(Y\) admits the representation \eqref{eq:spectral-representation-Y} through an orthogonal random measure \(Z_Y\), whose control measure is \(\mu_Y\). Motivated by the formal testing procedure in \eqref{YT-integral}, we define, for \(f\in L^2(\mathbb R)\),
$$
\widetilde Y(f)
:=
\int_{\mathbb R}\widehat f(\xi)\,Z_Y(\dd\xi),
$$
which defines the bounded linear operator \(\widetilde Y:L^2(\mathbb R)\longrightarrow L^2(\Omega)\),
called the bounded linear operator associated with \(Y\). This spectral construction does not require the time-domain integral in \eqref{YT-integral} to exist; its formal relation with the time-domain representation is given in \eqref{Y-MOTI}. The precise definition and properties of \(\widetilde Y\) are given in Definition~\ref{def:Ytilde}.

The extra boundedness of the spectral density plays an essential role in our setting and is also natural from the Fourier-domain point of view. Since multiplication by the spectral density \(S_Y\) defines a bounded operator on \(L^2(\mathbb R)\), the second-order structure induced by such process is compatible with the \(L^2(\mathbb R)\) test-function setting. This condition includes many familiar short-memory stationary models \cite{GinovyanSahakyan2017}. For example, if the covariance function is absolutely integrable, then the corresponding spectral density is bounded and continuous \cite{Lindgren2012}. On the other hand, stationary processes with unbounded spectral densities, and those whose spectral measures are not absolutely continuous, are not covered by the present framework. The boundedness of \(S_Y\) ensures that the stochastic integral defining \(\widetilde Y(f)\) is well defined for every \(f\in L^2(\mathbb R)\), and further implies that \(\widetilde Y\) is linear and continuous. Thus, rather than treating \(Y\) as an \(L^2(\mathbb R)\)-valued random element in the Bochner space \(L^2(\Omega,L^2(\mathbb R))\), we describe it through the operator \(\widetilde Y:L^2(\mathbb R)\to L^2(\Omega)\). The connection between this operator and the original stationary process \(Y\) will be made precise in the next section.


This viewpoint is motivated by the general theory of generalized stochastic processes, in particular It\^o's theory of stationary random distributions \cite{Ito1954}, where random objects are described through their action on deterministic test functions. Abstractly, such an object may be represented by a continuous linear mapping
$
\mathcal X:\mathcal V\longrightarrow L^2(\Omega),
$
with \(\mathcal V\) a deterministic test-function space. Here we take \(\mathcal V=L^2(\mathbb R)\), so that \(\mathcal X\) is precisely the  operator \(\widetilde Y\). Related formulations on locally compact Abelian groups use the Segal algebra \(S_0(G)\) as the test space \cite{FeichtingerHormann}. A related motivation comes from the treatment of Gaussian white noise in
\cite{YangChenQianZhao2025}, where harmonic and heat-flow representations
yield random solutions in suitable Bergman-type stochastic spaces, allowing
sparse expansions to be constructed by SPOAFD. In contrast, the present work
starts from a general second-order stationary process in place of Gaussian
white noise.

Our method is based on SPOAFD. We work in weighted Bergman spaces, where normalized reproducing kernels provide the dictionary for the adaptive decomposition. 
Deterministic POAFD in Bergman-type spaces has been studied in
\cite{QuDang2018,QuDang2019}. Qu and Dang \cite{QuDang2018} consider
pre-orthogonal approximation in weighted Hardy spaces, including Bergman
spaces on the unit disc and the standard Bergman space on the upper
half-plane, while \cite{QuDang2019} develops POAFD algorithms for weighted
Bergman spaces on the unit disc. In contrast, the present work treats
weighted Bergman spaces on both the upper and lower half-planes. For
\(\sigma\in\{+,-\}\) and \(z=x+iy\in\mathbb C^\sigma\), we consider
\[B_\alpha^2(\mathbb C^\sigma)
=
\left\{
F_\sigma\in\operatorname{Hol}(\mathbb C^\sigma):
\int_{\mathbb C^\sigma}
|F_\sigma(z)|^2
\left(
\frac{|y|}
{(1+|y|)^2+x^2}
\right)^\alpha
\,dA(z)
<\infty
\right\}.\]

Since \(\widetilde Y\) is defined on \(L^2(\mathbb R)\), we use the Cauchy kernels
$$
p_z(t)=\frac{1}{t-z},
\qquad z\in\mathbb C\setminus\mathbb R,
$$
as test functions. The upper and lower-half-plane kernels have Fourier transforms supported on the negative and positive frequency half-axes, respectively \cite{DangQianYou2011}, and together have dense linear span in \(L^2(\mathbb R)\). We therefore define
\[
U_+(\omega,z)=\widetilde Y(p_z),\quad z\in\mathbb C^+,\qquad
U_-(\omega,z)=\widetilde Y(p_z),\quad z\in\mathbb C^-.
\]
respectively, define random analytic functions on the two half-planes. 
Since, for any \(z\in\mathbb C\setminus\mathbb R\),
\[
\|p_z\|_{L^2(\mathbb R)}^2
=
\frac{\pi}{|\operatorname{Im}z|},
\]
the boundedness of \(\widetilde Y\) provides a second-moment estimate for
\(U_\pm(\omega,z)=\widetilde Y(p_z)\) of order
\(1/|\operatorname{Im}z|\). Therefore, the weight $\rho_\alpha(z)$ should compensate for this growth near the real axis and also provide sufficient decay at infinity. 
For \(\alpha>1\), the weight factor satisfies
\[
\int_{\mathbb C^\pm}
\frac{\rho_\alpha(z)}
{|\operatorname{Im}z|}
\,dA(z)
<\infty,\quad \rho_\alpha(z)
=
\left(
\frac{|y|}
{(1+|y|)^2+x^2}
\right)^\alpha,
\quad z=x+iy.
\]
Hence, the growth of the second moment of \(U_\pm(\omega,z)\) near the real axis is compensated by the weight, while the denominator in \(\rho_\alpha\) provides sufficient decay as \(|z|\to\infty\). The integrability, together with the analyticity of the sample functions, allows us to prove that \(U_+\) and \(U_-\) admit versions analytic on \(\mathbb C^+\) and \(\mathbb C^-\), respectively, for almost every \(\omega\in\Omega\), and that, for \(\alpha>1\),
\[
U_+\in L^2\bigl(\Omega,B_\alpha^2(\mathbb C^+)\bigr),
\qquad
U_-\in L^2\bigl(\Omega,B_\alpha^2(\mathbb C^-)\bigr).
\]
For \(\sigma\in\{+,-\}\), 
the corresponding Bochner space is
\begin{equation}\label{Bochner-space-sigma}
	\begin{aligned}
		L^2\bigl(\Omega,B_\alpha^2(\mathbb C^\sigma)\bigr)
		=
		\Bigl\{
		&U_\sigma:\Omega\to B_\alpha^2(\mathbb C^\sigma):
		U_\sigma\ \text{is a }
		B_\alpha^2(\mathbb C^\sigma)\text{-valued random variable, and}
		\\
		&\mathbb E
		\int_{\mathbb C^\sigma}
		|U_\sigma(\omega,z)|^2
		\left(
		\frac{|y|}
		{(1+|y|)^2+x^2}
		\right)^\alpha
		\,dA(z)
		<\infty
		\Bigr\},
	\end{aligned}
\end{equation}
where \(z=x+iy\).
Its norm is
\begin{equation}\label{Bochner-space-sigma-norm}
	\|U_\sigma\|_{L^2(\Omega,B_\alpha^2(\mathbb C^\sigma))}^2
	=
	\mathbb E
	\|U_\sigma\|_{B_\alpha^2(\mathbb C^\sigma)}^2
	=
	\mathbb E
	\int_{\mathbb C^\sigma}
	|U_\sigma(\omega,z)|^2
	\left(
	\frac{|y|}
	{(1+|y|)^2+x^2}
	\right)^\alpha
	\,dA(z).
\end{equation}

We next introduce, for \(\alpha>1\), the Hilbert direct sum
$
\mathcal B_\alpha
=
B_\alpha^2(\mathbb C^+)
\oplus
B_\alpha^2(\mathbb C^-),
$
identified with the space of functions on
\(\mathbb C^+\cup\mathbb C^-\) of the form
\begin{equation}\label{eq-F}
	F(z)
	=
	\begin{cases}
		F_+(z), & z\in\mathbb C^+,\\[2pt]
		F_-(z), & z\in\mathbb C^-,
	\end{cases}	
\end{equation}
where
\(F_+\in B_\alpha^2(\mathbb C^+)\) and
\(F_-\in B_\alpha^2(\mathbb C^-)\).
Its norm is
$
\|F\|_{\mathcal B_\alpha}^2
=
\|F_+\|_{B_\alpha^2(\mathbb C^+)}^2
+
\|F_-\|_{B_\alpha^2(\mathbb C^-)}^2.
$
We then combine the two random analytic functions as
\begin{equation}\label{Cauchy transform}
	U(\omega,z)
	=
	\begin{cases}
		U_+(\omega,z), & z\in\mathbb C^+,\\[2pt]
		U_-(\omega,z), & z\in\mathbb C^-.
	\end{cases}
\end{equation}
The associated Bochner space is
\[
L^2(\Omega,\mathcal B_\alpha)
=
\left\{
U:\Omega\to\mathcal B_\alpha:
U\ \text{is a }\mathcal B_\alpha\text{-valued random variable and}\
\mathbb E\|U\|_{\mathcal B_\alpha}^2<\infty
\right\},
\]
with norm
$
\|U\|_{L^2(\Omega,\mathcal B_\alpha)}^2
=
\mathbb E\|U\|_{\mathcal B_\alpha}^2.
$
We prove that
$U\in L^2(\Omega,\mathcal B_\alpha)$.
Thus, the object to be approximated is the random
analytic function $U$ in $L^2(\Omega,\mathcal B_\alpha)$.

We then apply SPOAFD independently in the upper and lower
half-planes, using the normalized reproducing kernels of
$B_\alpha^2(\mathbb C^+)$ and $B_\alpha^2(\mathbb C^-)$
as dictionary elements and establishing the boundary
vanishing condition in both half-planes.
The resulting orthonormal systems
$\{B_j^+\}_{j\geq1}$ and $\{B_\ell^-\}_{\ell\geq1}$
provide the SPOAFD approximation
\[
U_N(\omega,z)
=
\begin{cases}
	U_{N,+}(\omega,z), & z\in\mathbb C^+,\\[2pt]
	U_{N,-}(\omega,z), & z\in\mathbb C^-,
\end{cases}
\]
where
\[
\begin{aligned}
	U_{N,+}(\omega,\cdot)
	&=
	\sum_{j=1}^{N}
	\left\langle
	U(\omega,\cdot),B_j^+
	\right\rangle_{\mathcal B_\alpha}B_j^+,\\
	U_{N,-}(\omega,\cdot)
	&=
	\sum_{\ell=1}^{N}
	\left\langle
	U(\omega,\cdot),B_\ell^-
	\right\rangle_{\mathcal B_\alpha}B_\ell^-.
\end{aligned}
\]
and $B_j^+$ and $B_\ell^-$ are the orthonormal
functions defined in \eqref{eq:orthonormal-cross}.
The approximations converge in their respective
weighted Bergman spaces in mean square, which gives
\[
\lim_{N\to\infty}
\mathbb E\|U_N-U\|_{\mathcal B_\alpha}^2=0.
\]


Furthermore, the difference between the Cauchy
kernels on the upper and lower half-planes
gives the Poisson kernel
\[
P_y(t)=\frac{1}{2\pi i}
\left(\frac{1}{t-iy}-\frac{1}{t+iy}\right)
=\frac{y}{\pi(t^2+y^2)}.
\]
This relation gives
\[
V_y(\omega,t)
=
\frac{
	U_+(\omega,t+iy)-U_-(\omega,t-iy)
}{2\pi i}
=
\widetilde Y(P_y(t-\cdot))(\omega).
\]
Replacing $U_\pm$ by their SPOAFD approximations
$U_{N,\pm}$, we define
\[
V_{N,y}(\omega,t)
=
\frac{
	U_{N,+}(\omega,t+iy)
	-
	U_{N,-}(\omega,t-iy)
}{2\pi i}.
\]
We then establish the weak mean-square recovery
of the original stationary process:
for every $\gamma\in C_c^\infty(\mathbb R)$,
\[
\lim_{y\to0^+}\lim_{N\to\infty}
\mathbb E\left|
\int_{\mathbb R}\gamma(t)
\left[
V_{N,y}(\omega,t)-Y_t(\omega)
\right]\,dt
\right|^2=0.
\]

\medskip

\begin{minipage}{0.96\textwidth}
	\raggedright
	\setlength{\parindent}{0pt}
	\setlength{\parskip}{0pt}
	\setlength{\abovedisplayskip}{3pt}
	\setlength{\belowdisplayskip}{3pt}
	\setlength{\abovedisplayshortskip}{3pt}
	\setlength{\belowdisplayshortskip}{3pt}
	Our approach is divided into the following steps.
	
	\textbf{Whole-line stationary process.}
	\[
	Y=\{Y_t:t\in\mathbb R\},
	\qquad
	Y\notin L^2\bigl(\Omega,L^2(\mathbb R)\bigr)
	\quad\text{in general}.
	\]
	Thus, the whole sample path cannot in general be treated as a
	finite-energy random element of \(L^2(\mathbb R)\).
	
	\par\noindent\hfill{\Large\(\Downarrow\)}\hfill\mbox{}\par
	
	\textbf{Bounded linear operator associated with the stationary process Y.}
	Using the spectral representation of \(Y\), we define
	\[
	\widetilde Y:L^2(\mathbb R)\longrightarrow L^2(\Omega),
	\qquad
	\widetilde Y(f)
	=
	\int_{\mathbb R}\widehat f(\xi)\,Z_Y(\dd\xi).
	\]
	The bounded spectral density ensures that \(\widetilde Y\) is a
	well-defined bounded linear operator.
	
	\par\noindent\hfill{\Large\(\Downarrow\)}\hfill\mbox{}\par
	
	\textbf{Cauchy testing.}
	We test \(\widetilde Y\) with the Cauchy kernels
	\[
	p_z(t)=\frac{1}{t-z},
	\qquad
	z\in\mathbb C\setminus\mathbb R.
	\]
	The upper- and lower-half-plane kernels correspond to the negative
	and positive frequency half-axes, respectively, and together form
	a complete testing family in \(L^2(\mathbb R)\).
	
	\par\noindent\hfill{\Large\(\Downarrow\)}\hfill\mbox{}\par
	
	\textbf{Random analytic functions.}
	Define
	\[
	U_+(\omega,z)=\widetilde Y(p_z),
	~ z\in\mathbb C^+,
	~
	U_-(\omega,z)=\widetilde Y(p_z),
	~ z\in\mathbb C^-.
	\]
	For almost every \(\omega\in\Omega\),
	\(U_+(\omega,\cdot)\) and \(U_-(\omega,\cdot)\) are analytic on
	\(\mathbb C^+\) and \(\mathbb C^-\), respectively. They are combined as \eqref{Cauchy transform}
	
	\par\noindent\hfill{\Large\(\Downarrow\)}\hfill\mbox{}\par
	
	\textbf{Weighted Bergman--Bochner setting.}
	For \(\alpha>1\),
	\(	\mathcal B_\alpha
	=
	B_\alpha^2(\mathbb C^+)
	\oplus
	B_\alpha^2(\mathbb C^-).\)
	We prove that \[ U_+\in L^2\bigl(\Omega,B_\alpha^2(\mathbb C^+)\bigr),
	~
	U_-\in L^2\bigl(\Omega,B_\alpha^2(\mathbb C^-)\bigr),\]
	and hence \( U\in L^2(\Omega,\mathcal B_\alpha). \)
	Thus, \( U\), rather than the original sample path \(Y_t(\omega)\),
	is the object to be approximated.
	
	\par\noindent\hfill{\Large\(\Downarrow\)}\hfill\mbox{}\par
	
	\textbf{SPOAFD sparse representation.}
	SPOAFD is then applied in
	\(L^2(\Omega,\mathcal B_\alpha)\) using normalized weighted Bergman
	reproducing kernels to obtain an adaptive sparse expansion of
	\(U\).
	
	\par\noindent\hfill{\Large\(\Downarrow\)}\hfill\mbox{}\par
	
	\textbf{Weak mean-square convergence.}
	Using the Poisson kernel
	$P_y(t)=\frac{y}{\pi(t^2+y^2)}$
	and the SPOAFD expansions on both half-planes,
	we define $V_{N,y}$ as in
	\eqref{eq:SPOAFD-Poisson-approximation}. As $N\to\infty$ and subsequently
	$y\to0^+$,
	these approximations recover the original
	stationary process in the weak mean-square sense,
	as shown in \eqref{eq:SPOAFD-weak-recovery}.
	
\end{minipage}

\section{Stationary Processes and Their Associated Linear Operators}
\label{sec:generalized-random-operator}

To describe a stationary process on the whole real line in an
\(L^2(\mathbb R)\) framework, we consider its action on deterministic
functions rather than treating its sample paths as elements of
\(L^2(\mathbb R)\). Starting directly from the spectral representation
of the process, we associate it with  a bounded linear mapping from
\(L^2(\mathbb R)\) into \(L^2(\Omega)\).

We use the Fourier transform
\[
\widehat f(\xi)
=
\mathcal F f(\xi)
=
\frac{1}{\sqrt{2\pi}}
\int_{\mathbb R}
f(t)e^{-it\xi}\,\dd t,\qquad f\in L^2(\mathbb R),
\]
understood as the Plancherel extension \cite{SteinWeiss1971}, i.e., the \(L^2\)-limit of \(\mathcal{F}f_k\) for any sequence \(f_k \in L^1 \cap L^2\) converging to \(f\) in \(L^2\), which is independent of the choice of the sequence.


\begin{definition}
	\label{def:bounded-spectral-density-stationary-process}
	A jointly measurable, centered second-order stationary process
	\(Y=\{Y_t:t\in\mathbb R\}\)
	is called a \emph{stationary process with bounded spectral density} if its
	spectral measure is absolutely continuous with respect to Lebesgue measure
	and has a bounded spectral density
	\[
	S_Y\in L^\infty(\mathbb R),
	\qquad
	S_Y\geq0
	\quad \text{a.e.}
	\]
	Its spectral measure is then given by
	\begin{equation}\label{eq:spectral-measure}
		\mu_Y(\dd\xi)
		=
		S_Y(\xi)\,\dd\xi.
	\end{equation}
\end{definition}
The covariance function of \(Y\) is
\begin{equation}
	\E[Y_t\overline{Y_s}]
	=
	r_Y(t-s),
\end{equation}
and
\begin{equation}\label{FT-RY}
	r_Y(\tau)
	=
	\frac{1}{2\pi}
	\int_{\mathbb R}
	e^{-i\tau\xi}S_Y(\xi)\,\dd\xi.
\end{equation}
By the spectral representation theorem, also known as the Cram\'{e}r--Khinchin decomposition \cite{Bremaud2014}, there exists an orthogonal random measure \(Z_Y\) with  \(\mu_Y\) as control
measure such that
\begin{equation}\label{eq:spectral-representation-Y}
	Y_t
	=
	\frac{1}{\sqrt{2\pi}}
	\int_{\mathbb R}
	e^{-it\xi}\,Z_Y(\dd\xi),
	\qquad
	t\in\mathbb R,
\end{equation}
in the usual second-order sense.


Motivated by It\^{o}'s formulation of random distributions as continuous
linear mappings with values in a Hilbert space of second-order random
variables \cite{Ito1954}, we introduce the following operator \( \widetilde Y\)
associated  with \(Y\).

\begin{definition}\label{def:Ytilde}
	Let \(Y\) be the stationary process with bounded spectral density introduced in
	Definition~\ref{def:bounded-spectral-density-stationary-process}. For
	\(f\in L^2(\mathbb R)\), define
	\begin{equation}\label{DEF-YTILDE}
		\widetilde Y(f)
		:=
		\int_{\mathbb R}
		\widehat f(\xi)\,Z_Y(\mathrm d\xi),
	\end{equation}
	where \(Z_Y(\mathrm d\xi)\) is the orthogonal random measure appearing in the
	spectral representation \eqref{eq:spectral-representation-Y}.
	The mapping
	\[
	\widetilde Y:
	L^2(\mathbb R)
	\longrightarrow
	L^2(\Omega),
	\qquad
	f\longmapsto \widetilde Y(f),
	\]
	is called the linear operator induced by \(Y\).
\end{definition}

The notation \(\widetilde Y\) emphasizes that this operator is induced
by the original stationary process \(Y\), rather than introduced as an
independent stochastic object. Its definition in
\eqref{DEF-YTILDE} is motivated by the  idea of testing \(Y\)
against a deterministic  function \(f\) with finite-energy. Importantly,
this motivation does not require well defined of the pathwise integral
\begin{equation}\label{YT-integral}
	\int_{\mathbb R} f(t)Y_t\,\dd t,
\end{equation}
the integral in \eqref{YT-integral} is not used
to define \(\widetilde Y(f)\).
Indeed, if one were in a situation where the time-domain integral in
\eqref{YT-integral} and the interchange of integrations were justified,
then the spectral representation of \(Y\) would formally give
\begin{equation}\label{Y-MOTI}
	\begin{aligned}
		\int_{\mathbb R}f(t)Y_t\,\dd t
		&=
		\int_{\mathbb R}
		f(t)
		\left(
		\frac{1}{\sqrt{2\pi}}
		\int_{\mathbb R}
		e^{-it\xi}\,Z_Y(\dd\xi)
		\right)
		\dd t\\
		&=
		\int_{\mathbb R}
		\left(
		\frac{1}{\sqrt{2\pi}}
		\int_{\mathbb R}
		f(t)e^{-it\xi}\,\dd t
		\right)
		Z_Y(\dd\xi)
		=
		\int_{\mathbb R}
		\widehat f(\xi)\,Z_Y(\dd\xi)
		=
		\widetilde Y(f).
	\end{aligned}
\end{equation}
Thus, the spectral-domain definition \eqref{DEF-YTILDE} directly gives the rigorous way to interpret the formal testing procedure in \eqref{Y-MOTI}, without assuming that the sample-path integral \eqref{YT-integral} exists.

For a stationary process on the whole real line, global
\(L^2\)-integrability of its sample paths should not in general be
expected. As noted by Br\'emaud \cite[Sec.~3.4.1]{Bremaud2014},
almost surely, a trajectory of a stationary stochastic process belongs
to neither \(L^1(\mathbb R)\) nor \(L^2(\mathbb R)\), unless it is
identically null. This also emphasizes that stationarity and
sample-path integrability concern different aspects of a stochastic
process: In the present second-order setting, stationarity is expressed
through the translation-invariant mean and covariance structure,
whereas the condition \(Y(\omega,\cdot)\in L^2(\mathbb R) \)
is a global pathwise integrability requirement. In our treatment
\(L^2(\mathbb R)\) sample-path condition is not imposed here.

Relatively, the time-domain integral in \eqref{YT-integral} need not
exist pathwise. Throughout this paper, \(\widetilde Y(f)\) is defined
directly and rigorously by the spectral-domain expression
\eqref{DEF-YTILDE}; the formal calculation in \eqref{Y-MOTI} serves
only as motivation. Whenever the integral in \eqref{YT-integral}
is well defined and the interchange is justified,  \eqref{Y-MOTI} becomes valid.

\begin{proposition}\label{prop:Ytilde-bounded}
	The operator defined in \eqref{DEF-YTILDE} is centered, linear, and bounded. More precisely,
	\begin{equation}\label{eq:Ytilde-bound}
		\|\widetilde Y(f)\|_{L^2(\Omega)}
		\leq
		\|S_Y\|_{L^\infty(\mathbb R)}^{1/2}
		\|f\|_{L^2(\mathbb R)},
		\qquad
		f\in L^2(\mathbb R).
	\end{equation}
\end{proposition}


The operator $\widetilde Y$ introduced in Definition~\ref{def:Ytilde}
can be understood through the spectral-domain stochastic integral operator.
This relation also explains why $\widetilde Y$ is defined
on the fixed space $L^2(\mathbb R)$.
Let $\mu_Y$ be the spectral measure defined in
\eqref{eq:spectral-measure}, and define
\[
L^2(\mu_Y)
=
\left\{
\varphi:\mathbb R\to\mathbb C:
\int_{\mathbb R}
|\varphi(\xi)|^2\,\mu_Y(d\xi)<\infty
\right\},
\]
with norm \(\|\varphi\|_{L^2(\mu_Y)}^2
=
\int_{\mathbb R}
|\varphi(\xi)|^2\,\mu_Y(d\xi). \)
The stochastic integral operator associated with
$Z_Y$ is defined by
\begin{equation}\label{DEF-I}
	I:L^2(\mu_Y)\longrightarrow L^2(\Omega),
	\qquad
	I(\varphi)
	=
	\int_{\mathbb R}\varphi(\xi)\,Z_Y(d\xi).
\end{equation}
Such an operator is also considered in \cite{Ibragimov2019}.
By the isometry property of stochastic integration
with respect to $Z_Y$ in \eqref{eq:Z-isometry}, the operator $I$ satisfies
\begin{equation}\label{eq:I-isometry-norm}
	\|I(\varphi)\|_{L^2(\Omega)}^2
	=
	\|\varphi\|_{L^2(\mu_Y)}^2,
\end{equation}
so $I$ is a linear isometry.
Since $\mu_Y$ depends on $Y$, the domain $L^2(\mu_Y)$
of $I$ is process-dependent. In contrast,
the boundedness of $S_Y$ ensures that, for every
$f\in L^2(\mathbb R)$,
\[
\|\widehat f\|_{L^2(\mu_Y)}^2
=
\int_{\mathbb R}
|\widehat f(\xi)|^2S_Y(\xi)\,d\xi
\leq
\|S_Y\|_{L^\infty(\mathbb R)}
\|f\|_{L^2(\mathbb R)}^2.
\]
Thus $\widehat f\in L^2(\mu_Y)$, and
\begin{equation}\label{Y-I-COMPOSITION}
	\widetilde Y
	=
	I\circ\mathcal F.
\end{equation}
Hence $I$ acts on the process-dependent spectral space
$L^2(\mu_Y)$, whereas $\widetilde Y$ acts on the fixed
time-domain space $L^2(\mathbb R)$.



We next identify the covariance structure of \(\widetilde Y\).
For $f,g\in L^2(\mathbb R)$, the 
isometry \eqref{eq:Z-isometry} gives
\begin{equation}\label{eq:Ytilde-covariance}
	\E\bigl[
	\widetilde Y(f)
	\overline{\widetilde Y(g)}
	\bigr]
	=
	\int_{\mathbb R}
	\widehat f(\xi)
	\overline{\widehat g(\xi)}
	\,\mu_Y(\dd\xi)
	=
	\int_{\mathbb R}
	S_Y(\xi)
	\widehat f(\xi)
	\overline{\widehat g(\xi)}
	\,\dd\xi.
\end{equation}
Accordingly, define
\begin{equation}\label{eq:QY}
	Q_Y(f,g)
	:=
	\int_{\mathbb R}
	S_Y(\xi)
	\widehat f(\xi)
	\overline{\widehat g(\xi)}
	\,\dd\xi,
	\qquad
	f,g\in L^2(\mathbb R).
\end{equation}
Then \(\E\bigl[
\widetilde Y(f)
\overline{\widetilde Y(g)}
\bigr]
=
Q_Y(f,g).\)

\begin{proposition}\label{prop:QY}
	The form \(Q_Y\) defined in \eqref{eq:QY} is a bounded positive
	sesquilinear form on \(L^2(\mathbb R)\). More precisely,
	\[
	|Q_Y(f,g)|
	\leq
	\|S_Y\|_{L^\infty(\mathbb R)}
	\|f\|_{L^2(\mathbb R)}
	\|g\|_{L^2(\mathbb R)},
	\]
	and \(Q_Y(f,f)\geq0,
	~
	f,g\in L^2(\mathbb R).\)
\end{proposition}

The preceding construction is associated specifically with the
stationary process \(Y\). We next present a more general result. It shows that the covariance formulation above is not restricted to spectral covariance forms: every bounded positive sesquilinear form on \(L^2(\mathbb R)\) can be represented by a positive bounded covariance operator and realized as the covariance form of a centered bounded linear operator from \(L^2(\mathbb R)\) into \(L^2(\Omega)\).

\begin{proposition}\label{prop:abstract-covariance-realization}
	Let \(Q\) be a positive sesquilinear form on \(L^2(\mathbb R)\) such
	that there exists \(M>0\) satisfying
	\[
	|Q(f,g)|
	\leq
	M
	\|f\|_{L^2(\mathbb R)}
	\|g\|_{L^2(\mathbb R)},
	\qquad
	f,g\in L^2(\mathbb R).
	\]
	Then there exists a unique positive bounded self-adjoint operator \( 	C_Q:
	L^2(\mathbb R)
	\longrightarrow
	L^2(\mathbb R)\)
	such that
	\[
	Q(f,g)
	=
	\inner{C_Qf}{g}_{L^2(\mathbb R)},
	\qquad
	f,g\in L^2(\mathbb R).
	\]
	Moreover, \(Q\) can be realized as the covariance form of a centered
	bounded second-order random operator \(	\mathcal X_Q:
	L^2(\mathbb R)
	\longrightarrow
	L^2(\Omega). \)
\end{proposition}
We now return to the operator \(\widetilde Y\) induced by the
stationary process \(Y\). Since its covariance form \(Q_Y\) is bounded
and positive, the Riesz representation theorem guarantees the existence of a unique
positive bounded self-adjoint operator \(	C_{\widetilde Y}:
L^2(\mathbb R)
\longrightarrow
L^2(\mathbb R)\)
such that
\begin{equation}\label{eq:covariance-operator-Ytilde}
	\E\bigl[
	\widetilde Y(f)
	\overline{\widetilde Y(g)}
	\bigr]
	=
	Q_Y(f,g)
	=
	\inner{C_{\widetilde Y}f}{g}_{L^2(\mathbb R)},
	\qquad
	f,g\in L^2(\mathbb R).
\end{equation}
The operator \(C_{\widetilde Y}\) is called the covariance operator
of the   operator \(\widetilde Y\).
Its form can be identified explicitly in the Fourier domain. Let
\(M_{S_Y}\) denote the multiplication operator \((M_{S_Y}h)(\xi)
=
S_Y(\xi)h(\xi). \)
Then \(Q_Y(f,g)
=
\inner{M_{S_Y}\widehat f}
{\widehat g}_{L^2(\mathbb R)}. \)
Since the Fourier transform is unitary,
\begin{equation}\label{eq:covariance-operator-fourier}
	C_{\widetilde Y}
	=
	\mathcal F^{-1}
	M_{S_Y}
	\mathcal F.
\end{equation}
Equivalently,
\begin{equation}\label{eq:covariance-operator-multiplier}
	\widehat{C_{\widetilde Y}f}(\xi)
	=
	S_Y(\xi)\widehat f(\xi),
	\qquad
	f\in L^2(\mathbb R).
\end{equation}
The preceding discussion can be summarized as.
\begin{theorem}\label{thm:Y-induced-random-operator}
	Let \(Y=\{Y_t:t\in\mathbb R\}\) be a jointly measurable,
	centered second-order stationary process with spectral density
	\(S_Y\in L^\infty(\mathbb R)\), \(S_Y\geq0\) a.e.
	Let \(Z_Y\) be its orthogonal random measure with control measure
	\(\mu_Y(\dd\xi)=S_Y(\xi)\,\dd\xi\).
	Then \(Y\) induces the operator
	\[
	\widetilde Y:L^2(\mathbb R)\longrightarrow L^2(\Omega),
	\qquad
	\widetilde Y(f)
	=
	\int_{\mathbb R}\widehat f(\xi)\,Z_Y(\dd\xi).
	\]
	Its covariance form and covariance operator satisfy
	\[
	\E\bigl[\widetilde Y(f)\overline{\widetilde Y(g)}\bigr]
	=
	Q_Y(f,g)
	=
	\int_{\mathbb R}
	S_Y(\xi)\widehat f(\xi)
	\overline{\widehat g(\xi)}\,\dd\xi
	=
	\inner{C_{\widetilde Y}f}{g}_{L^2(\mathbb R)}.
	\]
	Moreover,
	\(C_{\widetilde Y}=\mathcal F^{-1}M_{S_Y}\mathcal F\),
	or equivalently,
	\(\widehat{C_{\widetilde Y}f}(\xi)
	=S_Y(\xi)\widehat f(\xi)\),
	for \(f\in L^2(\mathbb R)\).
\end{theorem}
\begin{remark}\label{Z}
	We recall the spectral measure, orthogonal random measure,
	and stochastic integral
	\cite{Lindgren2012,Ibragimov2019,LindquistPicci2007}.
	Let $\mu_Y$ be the spectral measure of
	$\{Y_t\}_{t\in\mathbb R}$, with distribution function
	$F_Y(\lambda)=\mu_Y((-\infty,\lambda])$.
	Since $\mu_Y(d\xi)=S_Y(\xi)\,d\xi$, we have
	$F_Y(\lambda)=\int_{-\infty}^{\lambda}S_Y(\xi)\,d\xi$.
	The orthogonal-random-measure and orthogonal-increment
	formulations are equivalent. An orthogonal-increment process
	$W$ defines $Z_Y$ by
	$Z_Y((a,b])=W(b)-W(a)$ for $a<b$,
	and conversely, $Z_Y$ determines $W$ up to an additive
	random variable \cite{LindquistPicci2007}.
	This formulation is also used for wide-sense stationary
	processes \cite{Ibragimov2019}.
	The orthogonal random measure $Z_Y$ in
	\eqref{eq:spectral-representation-Y} is countably additive
	in mean square and satisfies
	$\mathbb E[Z_Y(A)\overline{Z_Y(B)}]
	=\mu_Y(A\cap B)$
	for Borel sets $A,B\subset\mathbb R$.
	Thus, $Z_Y(A)$ and $Z_Y(B)$ are orthogonal in
	$L^2(\Omega)$ whenever $A\cap B=\varnothing$
	\cite{Ibragimov2019}.

	For $u,v\in L^2(\mu_Y)$, the isometry property
	of stochastic integration with respect to $Z_Y$
	\cite{Lindgren2012} can be expressed in terms of
	the operator $I$ defined in \eqref{DEF-I} as
	\begin{equation}\label{eq:Z-isometry}
		\mathbb E\left[I(u)\overline{I(v)}\right]
		=
		\int_{\mathbb R}
		u(\xi)\overline{v(\xi)}\,\mu_Y(d\xi)
		=
		\int_{\mathbb R}
		u(\xi)\overline{v(\xi)}S_Y(\xi)\,d\xi.
	\end{equation}

\end{remark}

\section{Analyticity of Random Cauchy Transforms on the Upper and Lower Half-Planes}

We use the bounded linear operator $\widetilde Y$ defined
in Definition~\ref{def:Ytilde} through \eqref{DEF-YTILDE}
to construct random Cauchy transforms on the upper and
lower half-planes. We first establish their second-moment
bounds and then prove their almost-sure analyticity.
These results provide a basis for studying the random
Cauchy transforms in weighted Bergman spaces.

For \(z\in\mathbb C\setminus\mathbb R\), let the Cauchy kernel $p_z$ be the test function
\[
p_z(t)=\frac{1}{t-z},
\qquad t\in\mathbb R.
\]
If $z=x+iy$, then
\[
\|p_z\|_{L^2(\R)}^2
=
\int_\R\frac{\dd t}{|t-z|^2}
=
\int_\R\frac{\dd t}{(t-x)^2+y^2}
=
\frac{\pi}{|y|}.
\]
Hence
$
p_z\in L^2(\R),
~ z\in\C\setminus\R.
$

\begin{definition}\label{def:cauchy-transform}
	Let \(\widetilde Y\) be the operator defined in Definition~\ref{def:Ytilde}.
	Define
	\[
	U_+(\omega,z)=\widetilde Y(p_z),
	~ z\in\mathbb C^+,
	\qquad
	U_-(\omega,z)=\widetilde Y(p_z),
	~ z\in\mathbb C^-.
	\]
	Define the random function \(U\) by \eqref{Cauchy transform}.
	Then \(U\) is called the random Cauchy transform associated with
	\(\widetilde Y\) on the upper and lower half-planes.
	
\end{definition}
We first derive a second-moment estimate for the
random Cauchy transforms using the boundedness
of the covariance operator $C_{\widetilde Y}$.
\begin{proposition}\label{second-order estimation}
	For every $z=x+iy\in\C\setminus\R$, \(  \E|\widetilde Y(p_z)|^2
	\le
	\frac{\pi\|C_{\widetilde Y}\|}{|y|}.\)
	Consequently,
	\[
	\E|U_+(x+iy)|^2\le \frac{\pi\|C_{\widetilde Y}\|}{y},
	~ y>0,
	~
	\E|U_-(x+iy)|^2\le \frac{\pi\|C_{\widetilde Y}\|}{|y|},
	~ y<0.
	\]
\end{proposition}

The preceding estimate controls the second moments
of the random Cauchy transforms. To establish their
almost-sure analyticity, we also need  the following lemma. 
\begin{lemma}\label{lem:weighted-path-integrability}\cite{Ito2006}
	If $\{Y_t\}_{t\in\mathbb R}$ is a jointly measurable
	second-order stationary process, then
	\[
	\mathbb E\left[
	\int_{\mathbb R}\frac{|Y_t|}{1+t^2}\,dt
	\right]
	\leq
	\pi\bigl(\mathbb E|Y_0|^2\bigr)^{1/2}
	<\infty.
	\]
	Consequently,
	\begin{equation}\label{Yt-condition}
		\int_{\mathbb R}
		\frac{|Y_t(\omega)|}{1+t^2}\,dt<\infty
		\qquad \text{a.s.}
	\end{equation}
\end{lemma}
The integrability condition \eqref{Yt-condition}
allows pathwise integration against differences of
Cauchy kernels and their derivatives. Together with
the  definition \ref{def:cauchy-transform}, it leads to the following analyticity result.

\begin{proposition}\label{prop:cauchy-analyticity}
	Let $Y=\{Y_t\}_{t\in\mathbb R}$ be the stationary
	process defined in Definition~\ref{def:bounded-spectral-density-stationary-process}.
	Let $\widetilde Y$ be the bounded linear operator
	defined in Definition~\ref{def:Ytilde}.
	With probability one, the functions
	$z\mapsto U_\pm(\omega,z)$
	are analytic on $\mathbb C^\pm$, respectively,
	and satisfy
	\[
	U_\pm'(\omega,z)
	=
	\int_{\mathbb R}
	\frac{Y_t(\omega)}{(t-z)^2}\,dt,
	\qquad z\in\mathbb C^\pm,
	\]
	where the integrals converge absolutely.
	Moreover, for each fixed $z\in\mathbb C^\pm$,
	\[
	U_\pm'(\cdot,z)
	=
	\widetilde Y\left(
	\frac{1}{(\cdot-z)^2}
	\right)
	\quad\text{in }L^2(\Omega).
	\]
\end{proposition}
Thus, the random Cauchy transforms admit almost-surely
analytic modifications on both half-planes and satisfy
the second-moment bounds established above.

\section{Completeness of Cauchy Testing Functions from Both Half-Planes}

We now examine the completeness of the Cauchy
testing functions used to construct the random
Cauchy transforms. Their Fourier transforms
have complementary frequency supports on
the upper and lower half-planes.
We use this property to establish the
completeness of the Cauchy kernel family
from both half-planes in $L^2(\mathbb R)$.

An equivalent identity to the following lemma
appears in \cite{DangQianYou2011}.
A proof based on contour integration and the residue theorem
is provided in the Appendix.

\begin{lemma}\label{FT-cauchy}
	Let $p_z(t)=(t-z)^{-1}$. If $z\in\C^+$, then
	\begin{equation}\label{eq:FT-cauchy-plus}
		\widehat{p_z}(\xi)
		=
		\sqrt{2\pi}\,i\,e^{-iz\xi}\1_{(-\infty,0)}(\xi).
	\end{equation}
	If $z\in\C^-$, then
	\begin{equation}\label{eq:FT-cauchy-minus}
		\widehat{p_z}(\xi)
		=
		-\sqrt{2\pi}\,i\,e^{-iz\xi}\1_{(0,\infty)}(\xi).
	\end{equation}
\end{lemma}

The preceding identities show that the Cauchy
kernels from the two half-planes have
complementary frequency supports.
This leads to the following completeness result.

\begin{theorem}\label{thm:two-sided-cauchy-complete}
	Let
	$
	p_z(t)=\frac{1}{t-z},~ z\in\mathbb C\setminus\mathbb R .
	$
	Then the upper and lower half-plane Cauchy family
	$
	\{p_z:z\in\mathbb C^+\cup\mathbb C^-\}
	$
	is complete in \(L^2(\mathbb R)\). More precisely,
	$
	\overline{\operatorname{span}}\{p_z:z\in\mathbb C^+\cup\mathbb C^-\}
	=
	L^2(\mathbb R).
	$
\end{theorem}

\begin{remark}
	The completeness of the Cauchy family also admits a
	boundary-value interpretation. For $f\in L^2(\mathbb R)$,
	define
	\[
	C^\pm f(z)
	=
	\frac{1}{2\pi i}
	\int_{\mathbb R}\frac{f(t)}{t-z}\,dt,
	\qquad z\in\mathbb C^\pm.
	\]
	By the Plemelj formulas \cite{Qian2015Book},
	the non-tangential boundary values satisfy
	\[
	C^\pm f(x\pm i0)
	=
	\pm\frac12 f(x)+\frac{i}{2}Hf(x).
	\]
	where $Hf$ denotes the Hilbert transform of $f$.
	Hence,
	$C^+f(x+i0)-C^-f(x-i0)=f(x)$.
	Thus, the Cauchy transforms on both half-planes
	recover $f$ from its boundary values,
	explaining why both half-planes are considered.
\end{remark}

The completeness of the Cauchy kernel family
from both half-planes shows that the values
of the bounded operator $\widetilde Y$
on these testing functions determine its
action on $L^2(\mathbb R)$.
Together with the analyticity and second-moment
bounds established earlier, this motivates
the study of random Cauchy transforms
in weighted Bergman spaces.

\section{Random Cauchy Transforms in Weighted Bergman Bochner Spaces}

By Proposition~\ref{prop:cauchy-analyticity},
the random Cauchy transforms constructed in the
preceding section admit almost-surely analytic
modifications on both half-planes.
Proposition~\ref{second-order estimation} also
provides their second-moment bounds.
We now introduce weighted Bergman spaces and
show that the random Cauchy transform belongs
to the corresponding weighted Bergman Bochner space.
We also establish the reproducing kernel properties
needed for the subsequent SPOAFD construction.

\subsection{Weighted Bergman Spaces on the Unit Disk}

We begin with the classical weighted Bergman space
on the unit disk, whose weight motivates the
construction of weighted Bergman spaces on
the upper and lower half-planes.
Let
$
\mathbb D=\{w\in\mathbb C: |w|<1\}
$
be the unit disk. For \(\alpha>-1\), we consider the weighted Bergman space
\cite{hedenmalm2012theory,Zhu2007} on \(\mathbb D\), with the normalization
\begin{equation}\label{standard weighted Bergman space}
	A^2_\alpha(\mathbb D)
	=
	\left\{
	f\in \operatorname{Hol}(\mathbb D):
	\int_{\mathbb D}|f(w)|^2(1-|w|^2)^\alpha\,dA(w)<\infty
	\right\}.
\end{equation}
The corresponding inner product is
\begin{equation}\label{standard weighted Bergman space-inner}
	\langle f,g\rangle_{A^2_\alpha(\mathbb D)}
	=
	\int_{\mathbb D}
	f(w)\overline{g(w)}(1-|w|^2)^\alpha\,dA(w),	
\end{equation}
and the induced norm is
\(\|f\|^2_{A^2_\alpha(\mathbb D)}
=
\int_{\mathbb D}|f(w)|^2(1-|w|^2)^\alpha\,dA(z).\)
The space \(A^2_\alpha(\mathbb D)\) is a reproducing kernel Hilbert space. Its reproducing kernel has the standard form
\begin{equation}
	K_b^{\mathbb D}(w)
	=
	\frac{\alpha+1}{\pi}
	\frac{1}{(1-w\overline b)^{\alpha+2}},
	\qquad z,b\in\mathbb D.
\end{equation}
Thus, for every \(f\in A^2_\alpha(\mathbb D)\) and every \(b\in\mathbb D\), one has the reproducing formula
\(f(b)
=
\langle f,K_b^{\mathbb D}(\cdot)\rangle_{A^2_\alpha(\mathbb D)}.\)

\subsection{Weighted Bergman Spaces on the Upper and Lower Half-Planes}

We now use the Cayley transform to obtain
suitable weights on the upper and lower
half-planes from the unit disk weight.
Let \(\mathbb D=\{w\in\mathbb C:|w|<1\},
~
\mathbb C^+=\{z\in\mathbb C:\operatorname{Im}z>0\}.\)
The Cayley transform from the upper half-plane onto the unit disk is
\[
w=\psi_+(z)=\frac{z-i}{z+i}.
\]
Following the standard definition of weighted Bergman spaces on the unit
disk \cite{hedenmalm2012theory,Zhu2007}, we consider the weight
\((1-|w|^2)^\alpha\), \(\alpha>-1\).
The normalization factor \(\alpha+1\) is omitted since it does not
affect the corresponding Bergman space.
For \(z=x+iy\in\mathbb C^+\), we have
\[
1-|w|^2
=
\frac{4\,\operatorname{Im}z}{|z+i|^2},
\qquad
(1-|w|^2)^\alpha
=
4^\alpha
\left(
\frac{\operatorname{Im}z}{|z+i|^2}
\right)^\alpha.
\]
Thus, up to a constant factor, the disk weight corresponds to
the following weight on the upper and lower half-planes:
\begin{equation}\label{half-plane weight}
	\rho_\alpha(z)
	=
	\begin{cases}
		\displaystyle
		\left(
		\frac{\operatorname{Im}z}{|z+i|^2}
		\right)^\alpha,
		& z\in\mathbb C^+,\\[8pt]
		\displaystyle
		\left(
		\frac{-\operatorname{Im}z}{|z-i|^2}
		\right)^\alpha,
		& z\in\mathbb C^-
	\end{cases}
	=
	\left(
	\frac{|y|}{(1+|y|)^2+x^2}
	\right)^\alpha,
	\qquad z=x+iy\in\mathbb C\setminus\mathbb R.
\end{equation}
For \(\sigma\in\{+,-\}\), let \(\mathbb C^\sigma\)
denote the corresponding half-plane.
The weighted Bergman space on \(\mathbb C^\sigma\) is defined by
\begin{equation}\label{weighted Bergman space-sigma}
	B_\alpha^2(\mathbb C^\sigma)
	=
	\left\{
	F_\sigma\in \operatorname{Hol}(\mathbb C^\sigma):
	\|F_\sigma\|_{B_\alpha^2(\mathbb C^\sigma)}<\infty
	\right\},
\end{equation}
where
\begin{equation}\label{weighted Bergman space-sigma-norm}
	\|F_\sigma\|_{B_\alpha^2(\mathbb C^\sigma)}^2
	=
	\int_{\mathbb C^\sigma}
	|F_\sigma(z)|^2\rho_\alpha(z)\,dA(z).
\end{equation}
Its inner product is
\[
\langle F_\sigma,G_\sigma\rangle_{B_\alpha^2(\mathbb C^\sigma)}
=
\int_{\mathbb C^\sigma}
F_\sigma(z)\overline{G_\sigma(z)}\rho_\alpha(z)\,dA(z).
\]
Equipped with this inner product, $B_\alpha^2(\mathbb C^\sigma)$
is complete and hence is a Hilbert space.

We identify
\begin{equation}\label{eq-B-alpha}
	\mathcal B_\alpha
	=
	B_\alpha^2(\mathbb C^+)\oplus B_\alpha^2(\mathbb C^-)
\end{equation}
with the space of functions on the disjoint union
\(\mathbb C^+\cup\mathbb C^-\). For \(F\in\mathcal B_\alpha\), write \(F\) as in \eqref{eq-F},
where \(F_\pm=F|_{\mathbb C^\pm}
\in B_\alpha^2(\mathbb C^\pm)\).
For \(F,G\in\mathcal B_\alpha\), define
\[
\langle F,G\rangle_{\mathcal B_\alpha}
=
\langle F_+,G_+\rangle_{B_\alpha^2(\mathbb C^+)}
+
\langle F_-,G_-\rangle_{B_\alpha^2(\mathbb C^-)}.
\]
Consequently,
\(\|F\|_{\mathcal B_\alpha}^2
=
\|F_+\|_{B_\alpha^2(\mathbb C^+)}^2
+
\|F_-\|_{B_\alpha^2(\mathbb C^-)}^2.\)

\begin{remark}\label{zero-extended}
	The spaces \(B_\alpha^2(\mathbb C^+)\) and
	\(B_\alpha^2(\mathbb C^-)\) are naturally identified isometrically with
	their zero-extended copies in \(\mathcal B_\alpha\). Define
	$$
	\iota_+(F_+)(z)
	=
	\begin{cases}
		F_+(z), & z\in\mathbb C^+,\\
		0, & z\in\mathbb C^-,
	\end{cases}
	\qquad
	\iota_-(F_-)(z)
	=
	\begin{cases}
		0, & z\in\mathbb C^+,\\
		F_-(z), & z\in\mathbb C^-.
	\end{cases}
	$$
	Then \(\iota_+\) and \(\iota_-\) are linear isometric isomorphisms onto
	\(B_\alpha^2(\mathbb C^+)\oplus\{0\}\) and
	\(\{0\}\oplus B_\alpha^2(\mathbb C^-)\), respectively, since
	\(
	\|\iota_+(F_+)\|_{\mathcal B_\alpha}
	=
	\|F_+\|_{B_\alpha^2(\mathbb C^+)},
	~
	\|\iota_-(F_-)\|_{\mathcal B_\alpha}
	=
	\|F_-\|_{B_\alpha^2(\mathbb C^-)}.
	\)
	Thus,
	\[\mathcal B_\alpha
	\cong
	\left(B_\alpha^2(\mathbb C^+)\oplus\{0\}\right)
	\oplus
	\left(\{0\}\oplus B_\alpha^2(\mathbb C^-)\right).\]
	Moreover, if \(F=\iota_+(F_+)\) and \(G=\iota_-(G_-)\), then
	\(
	\langle F,G\rangle_{\mathcal B_\alpha}
	=
	\langle F_+,0\rangle_{B_\alpha^2(\mathbb C^+)}
	+
	\langle 0,G_-\rangle_{B_\alpha^2(\mathbb C^-)}
	=0.
	\)
	Hence, the zero-extended upper and lower half-plane components are
	orthogonal in \(\mathcal B_\alpha\).
\end{remark}

\subsection{Bochner-Space Membership of the Random Cauchy Transforms $U$}

Having defined the weighted Bergman spaces
on both half-planes and their direct sum
$\mathcal B_\alpha$, we now show that the
random Cauchy transform $U$ defined in
\eqref{Cauchy transform} belongs to
$L^2(\Omega,\mathcal B_\alpha)$.
By Proposition~\ref{second-order estimation},
we have
\[
\mathbb E|U_\pm(z)|^2
\leq
\frac{\pi\|C_{\widetilde Y}\|}
{|\operatorname{Im}z|},
\qquad z\in\mathbb C^\pm.
\]
By Tonelli's theorem,
\[
\mathbb E\|U\|_{\mathcal B_\alpha}^2
=
\sum_{\sigma\in\{+,-\}}
\int_{\mathbb C^\sigma}
\mathbb E|U_\sigma(z)|^2
\rho_\alpha(z)\,\dd A(z)
\leq
\pi\|C_{\widetilde Y}\|
\sum_{\sigma\in\{+,-\}}
\int_{\mathbb C^\sigma}
\frac{\rho_\alpha(z)}
{|\operatorname{Im}z|}\,\dd A(z).
\]
Thus, it remains to establish the integrability of
$\rho_\alpha(z)/|\operatorname{Im}z|$ over both
half-planes, where $\rho_\alpha$ is defined in
\eqref{half-plane weight}.
The following lemma provides this result for $\alpha>1$.

\begin{lemma}\label{lem:weight-integrability}
	Let $\alpha>1$ and $\rho_\alpha$ be defined in
	\eqref{half-plane weight}. Then
	\[
	I_\alpha
	:=
	\int_{\mathbb C^+}
	\frac{\rho_\alpha(z)}{\operatorname{Im}z}
	\,dA(z)
	=
	c_\alpha B(\alpha,\alpha-1)
	<\infty,
	\]
	where
	$c_\alpha=\int_{\mathbb R}(1+u^2)^{-\alpha}\,du<\infty$
	and $B$ is the Beta function. And, \(	\int_{\mathbb C^-}
	\frac{\rho_\alpha(z)}{|\operatorname{Im}z|}
	\,dA(z)
	<\infty.\)
\end{lemma}

Combining Lemma~\ref{lem:weight-integrability}
with the second-moment estimate in
Proposition~\ref{second-order estimation},
we obtain the following  result.

\begin{theorem}\label{THM-Bochner}
	Let \(\widetilde Y\) be the bounded linear operator defined in
	Definition~\ref{def:Ytilde}, with covariance operator
	\(C_{\widetilde Y}\).
	If \(\alpha>1\), then
	\[
	U_+\in L^2\bigl(\Omega,B_\alpha^2(\C^+)\bigr),
	\qquad
	U_-\in L^2\bigl(\Omega,B_\alpha^2(\C^-)\bigr).
	\]
	Moreover,
	\[
	\E\|U_+\|_{B_\alpha^2(\C^+)}^2
	\leq
	\pi\|C_{\widetilde Y}\|I_\alpha,
	\qquad
	\E\|U_-\|_{B_\alpha^2(\C^-)}^2
	\leq
	\pi\|C_{\widetilde Y}\|I_\alpha.
	\]
	Let \(U\) be defined by \eqref{Cauchy transform}. Then \(U\in L^2(\Omega,\cB), \)	
	and
	\[
	\E\|U\|_{\cB}^2
	=
	\E\|U_+\|_{B_\alpha^2(\C^+)}^2
	+
	\E\|U_-\|_{B_\alpha^2(\C^-)}^2
	\leq
	2\pi\|C_{\widetilde Y}\|I_\alpha
	<
	\infty.
	\]
\end{theorem}

\subsection{Reproducing Kernels and Kernel Completeness}

Having established that
$U\in L^2(\Omega,\mathcal B_\alpha)$,
we now examine the reproducing kernel
structure of the weighted Bergman spaces.
In particular, we establish the completeness
of their reproducing kernel families, which
will be used to construct the normalized
kernel dictionary for SPOAFD.

\begin{lemma}\label{lem:bergman-rkhs}
	Let $\alpha>0$. For each $\sigma\in\{+,-\}$,
	the weighted Bergman space
	\[
	B^2_\alpha(\mathbb C^\sigma)
	=
	\left\{
	F_\sigma\in\operatorname{Hol}(\mathbb C^\sigma):
	\int_{\mathbb C^\sigma}
	|F_\sigma(z)|^2\rho_\alpha(z)\,dA(z)<\infty
	\right\}
	\]
	is a reproducing kernel Hilbert space.
	In particular, point evaluation is bounded, and
	for every $a\in\mathbb C^\sigma$, there exists a unique
	$K_a^{\mathbb C^\sigma}\in B^2_\alpha(\mathbb C^\sigma)$
	such that
	\[
	F_\sigma(a)
	=
	\left\langle F_\sigma,K_a^{\mathbb C^\sigma}\right\rangle_
	{B^2_\alpha(\mathbb C^\sigma)},
	\qquad F_\sigma\in B^2_\alpha(\mathbb C^\sigma).
	\]
	The reproducing kernel is given by
	$K^{\mathbb C^\sigma}(z,a)
	=K_a^{\mathbb C^\sigma}(z)$.
\end{lemma}

\begin{lemma}
	For each $\sigma\in\{+,-\}$,
	$
	\overline{\operatorname{span}}\{K^{\mathbb C^\sigma}_{a}~ on:a\in\C^\sigma\}
	=
	B^2_\alpha(\C^\sigma).
	$
	Equivalently, the normalized kernel dictionary
	\[
	\left\{\frac{K^{\mathbb C^\sigma}_{a}}{\|K^{\mathbb C^\sigma}_{a}\|_{B^2_\alpha(\C^\sigma)}}:a\in\C^\sigma\right\}
	\]
	is complete in $B^2_\alpha(\C^\sigma)$.
\end{lemma}

These results provide the Hilbert space setting
and complete reproducing kernel families needed
for SPOAFD.

\section{Normalized Bergman Kernel Dictionary and Boundary Vanishing Conditions}

To prepare for the SPOAFD construction, we first introduce
a normalized reproducing kernel dictionary in
$\mathcal B_\alpha$ and establish its completeness.
We then prove the boundary vanishing properties
in both $\mathcal B_\alpha$ and
$L^2(\Omega,\mathcal B_\alpha)$.

We consider \eqref{eq-B-alpha}
with the space of functions on the disjoint union
\(\mathbb C^+\cup\mathbb C^-\). Namely, every \(F\in\mathcal B_\alpha\) can be written
as in \eqref{eq-F}.
For \(a\in\mathbb C^+\), let \(K_a^{\mathbb C^+}\) denote the
reproducing kernel at \(a\) of the weighted Bergman space
\(B_\alpha^2(\mathbb C^+)\). For \(b\in\mathbb C^-\), let
\(K_b^{\mathbb C^-}\) denote the reproducing kernel at \(b\) of
\(B_\alpha^2(\mathbb C^-)\). Define
\[
e_a^+
=
\frac{K_a^{\mathbb C^+}}
{\left\|K_a^{\mathbb C^+}\right\|_{B_\alpha^2(\mathbb C^+)}},
~~
a\in\mathbb C^+,
\qquad
e_b^-
=
\frac{K_b^{\mathbb C^-}}
{\left\|K_b^{\mathbb C^-}\right\|_{B_\alpha^2(\mathbb C^-)}},
~~
b\in\mathbb C^-.
\]
We embed these normalized reproducing kernels into
\(\mathcal B_\alpha\) by setting
\begin{equation}\label{E+a}
	E_a^+(z)
	=
	\begin{cases}
		e_a^+(z), & z\in\mathbb C^+,\\[2pt]
		0, & z\in\mathbb C^-,
	\end{cases}
	~~
	a\in\mathbb C^+,
\end{equation}

\begin{equation}\label{E-b}
	E_b^-(z)
	=
	\begin{cases}
		0, & z\in\mathbb C^+,\\[2pt]
		e_b^-(z), & z\in\mathbb C^-,
	\end{cases}
	~~
	b\in\mathbb C^-.
\end{equation}
The normalized upper and lower half-plane Bergman kernel dictionary is
then defined by
\begin{equation}\label{eq:SPOAFD-dictionary}
	\mathcal D_\alpha
	=
	\mathcal D_\alpha^+
	\cup
	\mathcal D_\alpha^-
	=
	\left\{E_a^+:a\in\mathbb C^+\right\}
	\cup
	\left\{E_b^-:b\in\mathbb C^-\right\}	.
\end{equation}
By the completeness of the normalized reproducing kernel dictionaries
in the two half-plane Bergman spaces, \(\overline{\operatorname{span}}
\{e_a^+:a\in\mathbb C^+\}
=
B_\alpha^2(\mathbb C^+),
~
\overline{\operatorname{span}}
\{e_b^-:b\in\mathbb C^-\}
=
B_\alpha^2(\mathbb C^-).\)
Under the isometric identifications described in the Remark \ref{zero-extended},
the normalized kernels \(e_a^+\) and \(e_b^-\) are identified with their
zero-extended counterparts \(E_a^+\) and \(E_b^-\), respectively.
Hence, the completeness of the normalized kernel dictionaries in the two
half-plane Bergman spaces yields
\[
\overline{\operatorname{span}}
\{E_a^+:a\in\mathbb C^+\}
=
B_\alpha^2(\mathbb C^+)\oplus\{0\},
~
\overline{\operatorname{span}}
\{E_b^-:b\in\mathbb C^-\}
=
\{0\}\oplus B_\alpha^2(\mathbb C^-).
\]
It follows that \[\overline{\operatorname{span}}\mathcal D_\alpha
=
B_\alpha^2(\mathbb C^+)
\oplus
B_\alpha^2(\mathbb C^-)
=
\mathcal B_\alpha. \]

We next establish the boundary vanishing properties
of the normalized reproducing kernels, beginning
with the deterministic half-plane Bergman spaces.

\begin{lemma}[Boundary vanishing in the half-plane spaces]\label{lem:boundary-vanishing-half-plane-detailed}
	Let $\alpha>1$. For each $\sigma\in\{+,-\}$, let
	$K^{\mathbb C^\sigma}_{a}$ denote the reproducing kernel of $B^2_\alpha(\mathbb C^\sigma)$
	at the point $a\in\mathbb C^\sigma$, and set
	\[
	e_a^\sigma
	=
	\frac{K^{\mathbb C^\sigma}_{a}}{\|K^{\mathbb C^\sigma}_{a}\|_{B^2_\alpha(\C^\sigma)}}.
	\]
	Then for every $F_\sigma\in B^2_\alpha(\mathbb C^\sigma)$,
	\[
	\lim_{a\to \partial\mathbb C^\sigma}
	\left|
	\left\langle F_\sigma,e_a^\sigma\right\rangle_{B^2_\alpha(\mathbb C^\sigma)}
	\right|
	=0,
	\]
	where the boundary $\partial\mathbb C^\sigma$ is understood in the sense of the one-point compactification of $\mathbb C$, i.e., $\partial\mathbb C^\sigma = \mathbb R \cup \{\infty\}$.
\end{lemma}

\begin{corollary}\label{cor:boundary-vanishing}
	For every \(F\in\mathcal B_\alpha\), let
	\(F_+\) and \(F_-\) be defined as in \eqref{eq-F}.
	Then
	\[
	\lim_{a\to\partial\mathbb C^+}
	\left|
	\left\langle F,E_a^+\right\rangle_{\mathcal B_\alpha}
	\right|
	=0,
	\qquad
	\lim_{b\to\partial\mathbb C^-}
	\left|
	\left\langle F,E_b^-\right\rangle_{\mathcal B_\alpha}
	\right|
	=0.
	\]
\end{corollary}

We now extend the boundary vanishing property
to the weighted Bergman Bochner space
$L^2(\Omega,\mathcal B_\alpha)$.

\begin{lemma}\label{lem:random-boundary-vanishing}
	Let \( U\in L^2(\Omega,\mathcal B_\alpha)\). Then
	\[
	\lim_{a\to\partial\mathbb C^+}
	\mathbb E
	\left|
	\left\langle  U,E_a^+\right\rangle_{\mathcal B_\alpha}
	\right|^2
	=0,
	\qquad
	\lim_{b\to\partial\mathbb C^-}
	\mathbb E
	\left|
	\left\langle  U,E_b^-\right\rangle_{\mathcal B_\alpha}
	\right|^2
	=0.
	\]
\end{lemma}

These boundary vanishing properties will be
used in the next section to establish the
existence of maximizers for the SPOAFD
selection functionals.

\section{Stochastic Pre-Orthogonal Adaptive Fourier Decomposition in the
	Weighted Bergman Bochner Space}

Using the normalized kernel dictionary and boundary
vanishing properties established above, we now
construct SPOAFD in $L^2(\Omega,\mathcal B_\alpha)$,
following the approach in \cite{Qian2022}.
We apply SPOAFD independently to $U_+$ and $U_-$
in their respective weighted Bergman Bochner spaces.
The two approximations will then be combined in
$L^2(\Omega,\mathcal B_\alpha)$.

\medskip
\noindent\textbf{Step 1. Dictionary and previously selected elements}

Let $\mathcal D_\alpha$ be the normalized Bergman kernel
dictionary in $\mathcal B_\alpha$ defined in
\eqref{eq:SPOAFD-dictionary}. We use the two kernel families
\[
\mathcal D_\alpha^+=
\{E_a^+:a\in\mathbb C^+\},
\qquad
\mathcal D_\alpha^-=\{E_b^-:b\in\mathbb C^-\},
\]
separately. Suppose that, before the $n$-th selection
in the upper and lower half-planes, respectively,
\[
E_{a_1}^+,\ldots,E_{a_{n-1}}^+
\qquad\text{and}\qquad
E_{b_1}^-,\ldots,E_{b_{n-1}}^-
\]
have been selected. These two sequences are constructed
independently.

The following orthogonality property allows the two half-plane
components to be treated separately.

\begin{theorem}\label{thm:half-plane-orthogonality}
	The upper and lower half-plane kernel families generate orthogonal
	subspaces of $\mathcal B_\alpha$. In particular,
	\begin{equation}\label{eq:upper-lower-subspaces}
		\overline{\operatorname{span}}
		\{E_a^+:a\in\mathbb C^+\}
		=
		B_\alpha^2(\mathbb C^+)\oplus\{0\},
		\qquad
		\overline{\operatorname{span}}
		\{E_b^-:b\in\mathbb C^-\}
		=
		\{0\}\oplus B_\alpha^2(\mathbb C^-),
	\end{equation}
	and hence
	\begin{equation}\label{eq:upper-lower-orthogonality}
		\left\langle E_a^+,E_b^-\right\rangle_{\mathcal B_\alpha}=0,
		\qquad
		a\in\mathbb C^+,\quad b\in\mathbb C^-.
	\end{equation}
	Consequently, the Gram--Schmidt orthogonalization can be carried out
	separately in the two half-plane components.
\end{theorem}

\medskip
\noindent\textbf{Step 2. Orthogonalization of the previously selected elements}

Let $\{B_1^+,\ldots,B_{n-1}^+\}$ be the orthonormal system
obtained from $\{E_{a_1}^+,\ldots,E_{a_{n-1}}^+\}$,
and let $\{B_1^-,\ldots,B_{n-1}^-\}$ be the orthonormal system
obtained from $\{E_{b_1}^-,\ldots,E_{b_{n-1}}^-\}$.
Here the orthogonalization is performed independently in each
half-plane. By Theorem~\ref{thm:half-plane-orthogonality},
\begin{equation}\label{eq:orthonormal-cross}
	\left\langle B_j^+,B_\ell^-\right\rangle_{\mathcal B_\alpha}
	=0,
	\qquad
	1\leq j,\ell\leq n-1.
\end{equation}
Consequently, the two systems
\begin{equation}\label{eq:SPOAFD-orthonormal-system}
	B_1^+,\ldots,B_{n-1}^+,
	\qquad
	B_1^-,\ldots,B_{n-1}^-
\end{equation}
are orthogonal to each other in $\mathcal B_\alpha$.

\medskip

\noindent\textbf{Step 3. Pre-orthogonalization of a new candidate}

First, let $d=E_a^+$, $a\in\mathbb C^+$.
Denote by $\widetilde E_a^{+}$ the associated kernel or
multiple kernel at the $n$-th step of the upper half-plane
procedure. If $a$ has not been selected before, then
$\widetilde E_a^{+}=E_a^+$; otherwise, it is the
corresponding embedded multiple kernel obtained by
differentiation with respect to the conjugate parameter.

Since $\widetilde E_a^{+}$ and $B_j^+$ vanish in
$\mathbb C^-$, the orthogonalization is
\begin{equation}\label{eq:SPOAFD-Qplus}
	Q_n(E_a^+)
	=
	\widetilde E_a^{+}
	-
	\sum_{j=1}^{n-1}
	\left\langle
	\widetilde E_a^{+},B_j^+
	\right\rangle_{\mathcal B_\alpha}
	B_j^+.
\end{equation}
Similarly, for $d=E_b^-$, $b\in\mathbb C^-$,
let $\widetilde E_b^{-}$ denote the corresponding
kernel or multiple kernel in the lower half-plane.
The orthogonalization is
\begin{equation}\label{eq:SPOAFD-Qminus}
	Q_n(E_b^-)
	=
	\widetilde E_b^{-}
	-
	\sum_{\ell=1}^{n-1}
	\left\langle
	\widetilde E_b^{-},B_\ell^-
	\right\rangle_{\mathcal B_\alpha}
	B_\ell^-.
\end{equation}
For $a\in\mathbb C^+$ and $b\in\mathbb C^-$,
define the normalized pre-orthogonalized candidates by
\begin{equation}\label{eq:SPOAFD-candidate-B}
	B_n^+(a)
	=
	\frac{Q_n(E_a^+)}
	{\left\|Q_n(E_a^+)\right\|_{\mathcal B_\alpha}},
	\qquad
	B_n^-(b)
	=
	\frac{Q_n(E_b^-)}
	{\left\|Q_n(E_b^-)\right\|_{\mathcal B_\alpha}},
\end{equation}
whenever the corresponding denominators are nonzero.

Once the maximizing parameters $a_n$ and $b_n$
have been selected, the next orthonormal elements
are defined by
\begin{equation}\label{eq:SPOAFD-Bplus}
	B_n^+=B_n^+(a_n),
\end{equation}
and
\begin{equation}\label{eq:SPOAFD-Bminus}
	B_n^-=B_n^-(b_n).
\end{equation}
The selections of $a_n$ and $b_n$ are independent.

\medskip
\noindent\textbf{Step 4. Random residual}

Let $U\in L^2(\Omega,\mathcal B_\alpha)$.
Define the random residual by
\begin{equation}\label{eq:SPOAFD-residual}
	G_n(\omega,z)
	=
	\begin{cases}
		G_{n,+}(\omega,z), & z\in\mathbb C^+,\\[2pt]
		G_{n,-}(\omega,z), & z\in\mathbb C^-,
	\end{cases}
\end{equation}
where
\[
\begin{aligned}
	G_{n,+}(\omega,\cdot)
	&=
	U_+(\omega,\cdot)
	-
	\sum_{j=1}^{n-1}
	\left\langle
	U(\omega,\cdot),B_j^+
	\right\rangle_{\mathcal B_\alpha}
	B_j^+,\\
	G_{n,-}(\omega,\cdot)
	&=
	U_-(\omega,\cdot)
	-
	\sum_{\ell=1}^{n-1}
	\left\langle
	U(\omega,\cdot),B_\ell^-
	\right\rangle_{\mathcal B_\alpha}
	B_\ell^-.
\end{aligned}
\]
Thus,
$G_{n,+}\in L^2(\Omega,B_\alpha^2(\mathbb C^+))$
and
$G_{n,-}\in L^2(\Omega,B_\alpha^2(\mathbb C^-))$,
with $G_n\in L^2(\Omega,\mathcal B_\alpha)$.

\medskip
\noindent\textbf{Step 5. Selection functional}

By the direct-sum structure, each candidate interacts only with
the residual component in the same half-plane.
For $a\in\mathbb C^+$, orthogonalization and
normalization preserve the vanishing property
$B_n^+(a)(z)=0$ for $z\in\mathbb C^-$.
Similarly, for $b\in\mathbb C^-$,
$B_n^-(b)(z)=0$ for $z\in\mathbb C^+$.

We define the selection functionals for the
upper and lower half-planes, respectively, by
\begin{equation}\label{eq:SPOAFD-Jplus}
	J_n^+(a)
	=
	\mathbb E
	\left|
	\left\langle
	G_n(\omega,\cdot),B_n^+(a)
	\right\rangle_{\mathcal B_\alpha}
	\right|^2,
	\qquad a\in\mathbb C^+,
\end{equation}
and
\begin{equation}\label{eq:SPOAFD-Jminus}
	J_n^-(b)
	=
	\mathbb E
	\left|
	\left\langle
	G_n(\omega,\cdot),B_n^-(b)
	\right\rangle_{\mathcal B_\alpha}
	\right|^2,
	\qquad b\in\mathbb C^-.
\end{equation}
In the SPOAFD procedure, $J_n^+$ is maximized only over
$\mathbb C^+$, and $J_n^-$ is maximized only over
$\mathbb C^-$.

\medskip
\noindent\textbf{Step 6. Maximization and selection of the $n$-th element}

The existence of a maximizing dictionary element in each
half-plane follows from the corresponding boundary vanishing
property.

\begin{theorem}\label{thm:SPOAFD-selection-existence}
	By Corollary~\ref{cor:boundary-vanishing},
	Lemma~\ref{lem:random-boundary-vanishing}, and the continuity
	of $J_n^+$ and $J_n^-$, we have
	\begin{equation}\label{eq:SPOAFD-J-boundary}
		\lim_{a\to\partial\mathbb C^+}J_n^+(a)=0,
		\qquad
		\lim_{b\to\partial\mathbb C^-}J_n^-(b)=0.
	\end{equation}
	Hence both functions attain their maxima at
	\begin{equation}\label{eq:SPOAFD-maximizers}
		a_n
		\in
		\operatorname*{arg\,max}_{a\in\mathbb C^+}J_n^+(a),
		\qquad
		b_n
		\in
		\operatorname*{arg\,max}_{b\in\mathbb C^-}J_n^-(b).
	\end{equation}
	Moreover, the two residuals are considered independently.
	If
	\begin{equation}\label{eq:SPOAFD-nonzero-residual}
		G_{n,+}\neq0
		\quad\text{in }L^2(\Omega,B_\alpha^2(\mathbb C^+)),
	\end{equation}
	then
	\begin{equation}\label{eq:SPOAFD-positive-sup}
		\sup_{a\in\mathbb C^+}J_n^+(a)>0.
	\end{equation}
	Similarly, if
	$G_{n,-}\neq0$ in
	$L^2(\Omega,B_\alpha^2(\mathbb C^-))$,
	then $\sup_{b\in\mathbb C^-}J_n^-(b)>0$.
\end{theorem}

The maximizing parameters $a_n$ and $b_n$ in
\eqref{eq:SPOAFD-maximizers} determine $B_n^+$ and $B_n^-$
through \eqref{eq:SPOAFD-Bplus} and \eqref{eq:SPOAFD-Bminus},
respectively.

\subsection{Energy Identity and Convergence}

Let
$\{B_j^+\}_{j\geq1}$ and $\{B_\ell^-\}_{\ell\geq1}$
be the corresponding orthonormal systems obtained
from the maximal selection rules in
\eqref{eq:SPOAFD-maximizers}.
By \eqref{eq:orthonormal-cross},
these two systems are orthogonal in
$\mathcal B_\alpha$.

For $N\geq1$, define the SPOAFD approximation by
\begin{equation}\label{eq:SPOAFD-partial-sum-1}
	U_N(\omega,z)
	=
	\begin{cases}
		U_{N,+}(\omega,z), & z\in\mathbb C^+,\\[2pt]
		U_{N,-}(\omega,z), & z\in\mathbb C^-,
	\end{cases}
\end{equation}
where
\[
\begin{aligned}
	U_{N,+}(\omega,\cdot)
	&=
	\sum_{j=1}^{N}
	\left\langle
	U(\omega,\cdot),B_j^+
	\right\rangle_{\mathcal B_\alpha}
	B_j^+,\\
	U_{N,-}(\omega,\cdot)
	&=
	\sum_{\ell=1}^{N}
	\left\langle
	U(\omega,\cdot),B_\ell^-
	\right\rangle_{\mathcal B_\alpha}
	B_\ell^-.
\end{aligned}
\]
Thus, $U_N$ is formed by the $N$-term SPOAFD
approximations in the upper and lower half-planes,
and $U_N\in L^2(\Omega,\mathcal B_\alpha)$.

By the definition of the residual in
\eqref{eq:SPOAFD-residual},
\[
G_{N+1}(\omega,z)
=
\begin{cases}
	U_+(\omega,z)-U_{N,+}(\omega,z),
	& z\in\mathbb C^+,\\[2pt]
	U_-(\omega,z)-U_{N,-}(\omega,z),
	& z\in\mathbb C^-.
\end{cases}
\]
Therefore,
\begin{equation}\label{eq:G-U-U}
	G_{N+1}(\omega,z)
	=
	U(\omega,z)-U_N(\omega,z).
\end{equation}

\begin{theorem}\label{thm-spoafd-convergence}
	Let $U\in L^2(\Omega,\mathcal B_\alpha)$,
	and let $\{B_j^+\}_{j\geq1}$ and
	$\{B_\ell^-\}_{\ell\geq1}$
	be the orthonormal systems obtained from the
	SPOAFD procedures in the upper and lower
	half-planes, respectively.
	Let $U_{N,+}$, $U_{N,-}$, and $U_N$ be
	defined as in \eqref{eq:SPOAFD-partial-sum-1}.
	For every $N\geq1$, we have
	\[
	\mathbb E
	\|U_+-U_{N,+}\|_{B_\alpha^2(\mathbb C^+)}^2
	=
	\mathbb E
	\|U_+\|_{B_\alpha^2(\mathbb C^+)}^2
	-
	\sum_{j=1}^{N}
	\mathbb E
	\left|
	\left\langle U,B_j^+
	\right\rangle_{\mathcal B_\alpha}
	\right|^2,
	\]
	and
	\[
	\mathbb E
	\|U_--U_{N,-}\|_{B_\alpha^2(\mathbb C^-)}^2
	=
	\mathbb E
	\|U_-\|_{B_\alpha^2(\mathbb C^-)}^2
	-
	\sum_{\ell=1}^{N}
	\mathbb E
	\left|
	\left\langle U,B_\ell^-
	\right\rangle_{\mathcal B_\alpha}
	\right|^2.
	\]
	Moreover,
	\[
	\lim_{N\to\infty}
	\mathbb E
	\|U_+-U_{N,+}\|_{B_\alpha^2(\mathbb C^+)}^2
	=0
	\]
	and
	\[
	\lim_{N\to\infty}
	\mathbb E
	\|U_--U_{N,-}\|_{B_\alpha^2(\mathbb C^-)}^2
	=0.
	\]
	By the direct-sum structure of $\mathcal B_\alpha$,
	\[
	\mathbb E\|U-U_N\|_{\mathcal B_\alpha}^2
	=
	\mathbb E
	\|U_+-U_{N,+}\|_{B_\alpha^2(\mathbb C^+)}^2
	+
	\mathbb E
	\|U_--U_{N,-}\|_{B_\alpha^2(\mathbb C^-)}^2.
	\]
	Consequently,
	\[
	\lim_{N\to\infty}
	\mathbb E\|U-U_N\|_{\mathcal B_\alpha}^2
	=0.
	\]
	Equivalently,
	\begin{equation}\label{eq:U-expansion}
		U(\omega,z)
		=
		\begin{cases}
			\displaystyle
			\sum_{j=1}^{\infty}
			\left\langle
			U(\omega,\cdot),B_j^+
			\right\rangle_{\mathcal B_\alpha}
			B_j^+(z),
			& z\in\mathbb C^+,\\[12pt]
			\displaystyle
			\sum_{\ell=1}^{\infty}
			\left\langle
			U(\omega,\cdot),B_\ell^-
			\right\rangle_{\mathcal B_\alpha}
			B_\ell^-(z),
			& z\in\mathbb C^-,
		\end{cases}
	\end{equation}
	where the expansions converge in
	$L^2(\Omega,\mathcal B_\alpha)$.
\end{theorem}

\subsection{Coefficient Covariance Identity}

To connect the two SPOAFD expansions with the original operator \(\widetilde Y\), we identify their random coefficients. For \(\alpha>1\) and every \(F\in\mathcal B_\alpha\), we show that there exists a deterministic function \(\Phi_F\in L^2(\mathbb R)\), obtained from weighted integrals of the Cauchy kernels over the upper and lower half-planes, such that
\[
\left\langle U(\omega,\cdot),F(\cdot)\right\rangle_{\mathcal B_\alpha}
=
\widetilde Y(\Phi_F)(\omega)
\]
as an identity in \(L^2(\Omega)\).  
In particular, taking $F=B_j^+$ or $F=B_\ell^-$
gives the random coefficients in the corresponding
SPOAFD expansions.
For $F\in\cB$, write $F$ as in \eqref{eq-F},
where
$
F_+\in B_\alpha^2(\mathbb C^+),
~
F_-\in B_\alpha^2(\mathbb C^-).
$
Consider
\begin{equation}\label{Phi-F}
	\Phi_F^\pm(t)
	:=
	\int_{\mathbb C^\pm}
	\frac{\overline{F_\pm(z)}\rho_\alpha(z)}{t-z}\,\dd A(z),
	\qquad
	\Phi_F:=\Phi_F^++\Phi_F^-.
\end{equation}
The following lemma shows that these scalar integrals converge
absolutely for almost every $t\in\R$ and define elements of $L^2(\R)$.

\begin{lemma}\label{PHI-F}
	Let $\alpha>1$ and $F\in\cB$. Then $\Phi_F\in L^2(\R)$ with
	$\|\Phi_F\|_{L^2(\R)}
	\le\sqrt{2\pi I_\alpha}\,\|F\|_{\cB}$, where
	\[
	I_\alpha
	:=
	\int_{\Cp}
	\frac{\rho_\alpha(z)}{\operatorname{Im}z}\,\dd A(z)
	=
	\int_{\Cm}
	\frac{\rho_\alpha(z)}{|\operatorname{Im}z|}\,\dd A(z)
	<\infty.
	\]
	Moreover,
	\begin{equation}\label{inner-U-B}
		\inner{U}{F}_{\cB}={\widetilde Y}(\Phi_F)
	\end{equation}
	in $L^2(\Omega)$.
\end{lemma}

Using this representation, we obtain the second-moment identity for the random coefficients.

\begin{theorem}\label{thm:PHI-F}
	Let $\alpha>1$, let $F\in\cB$, and  let
	$\Phi_F$ be defined as in \eqref{Phi-F}. Then
	\[
	\E\left|\inner{U}{F}_{\cB}\right|^2
	=
	\inner{C_{\widetilde Y}\Phi_F}{\Phi_F}_{L^2(\R)}.
	\]
	Since
	$
	\widehat{C_{\widetilde Y}f}(\xi)
	=
	S_Y(\xi)\widehat f(\xi),
	$
	we further have
	\[
	\E\left|\inner{U}{F}_{\cB}\right|^2
	=
	\int_\R
	S_Y(\xi)|\widehat{\Phi_F}(\xi)|^2\,\dd\xi.
	\]
\end{theorem}
Thus, the second moments of the coefficients in both SPOAFD expansions are determined by the spectral density of the underlying stationary process.

\section{Weak Mean-Square Convergence to the Original Stationary Process}


We now show that the SPOAFD expansion of the random
Cauchy transforms provides a weak mean-square recovery
of the original stationary process.

Let
\[
P_y(t)=\frac{1}{\pi}\frac{y}{t^2+y^2},
\qquad y>0,
\]
be the Poisson kernel on $\mathbb R$.
The corresponding Poisson regularization can be
expressed through the random Cauchy transforms as
\begin{equation}\label{eq:Poisson-regularization}
	V_y(\omega,t)
	=
	\frac{
		U_+(\omega,t+iy)-U_-(\omega,t-iy)
	}{2\pi i}
	=
	\widetilde Y(P_y(t-\cdot))(\omega).
\end{equation}
where the last equality holds in $L^2(\Omega)$
for each fixed $t$.
Replacing $U_\pm$ by their SPOAFD approximations
$U_{N,\pm}$, we obtain
\begin{equation}\label{eq:SPOAFD-Poisson-approximation}
	V_{N,y}(\omega,t)
	=
	\frac{
		U_{N,+}(\omega,t+iy)
		-
		U_{N,-}(\omega,t-iy)
	}{2\pi i}.
\end{equation}
Using the mean-square convergence of SPOAFD,
the reproducing property of the weighted Bergman
spaces, and the approximate identity property
of the Poisson kernel, we establish the weak
mean-square recovery of the original stationary
process against test functions in
$C_c^\infty(\mathbb R)$.

\begin{theorem}\label{thm:SPOAFD-weak-recovery}
	Let $Y=\{Y_t\}_{t\in\mathbb R}$ be the stationary
	process defined in Definition~\ref{def:bounded-spectral-density-stationary-process}.
	Let $\widetilde Y$ be the bounded linear operator
	defined in Definition~\ref{def:Ytilde}.
	For $\alpha>1$, let
	$U\in L^2(\Omega,\mathcal B_\alpha)$
	be the random Cauchy transforms.
	Let
	\[
	U_N(\omega,z)
	=
	\begin{cases}
		U_{N,+}(\omega,z), & z\in\mathbb C^+,\\
		U_{N,-}(\omega,z), & z\in\mathbb C^-,
	\end{cases}
	\]
	where
	\[
	\begin{aligned}
		U_{N,+}(\omega,z)
		&=
		\sum_{j=1}^{N}
		\left\langle
		U(\omega,\cdot),B_j^+
		\right\rangle_{\mathcal B_\alpha}
		B_j^+(z),
		\qquad z\in\mathbb C^+,\\
		U_{N,-}(\omega,z)
		&=
		\sum_{\ell=1}^{N}
		\left\langle
		U(\omega,\cdot),B_\ell^-
		\right\rangle_{\mathcal B_\alpha}
		B_\ell^-(z),
		\qquad z\in\mathbb C^-.
	\end{aligned}
	\]
	Then \(U_N\) is the approximation obtained from the \(N\)-term SPOAFD expansions in the upper and lower half-planes, satisfying
	\[
	\lim_{N\to\infty}
	\mathbb E\|U_N-U\|_{\mathcal B_\alpha}^2=0.
	\]
	Moreover, for every $\gamma\in C_c^\infty(\mathbb R)$,
	\begin{equation}\label{eq:SPOAFD-weak-recovery}
		\lim_{y\to0^+}\lim_{N\to\infty}
		\mathbb E\Bigg|
		\int_{\mathbb R}\gamma(t)
		\Bigg[
		V_{N,y}(\omega,t)
		-Y_t(\omega)
		\Bigg]\,dt
		\Bigg|^2=0.
	\end{equation}
\end{theorem}

\section*{Acknowledgement}
The research is supported by National Science and Technology Major Project (2027ZD01999504), the Major Project of Guangzhou National Laboratory (Grant No. GZNL2024A01004), the Science and Technology Development Fund of Macau SAR (Grant No. 0015/2025/AFJ), and the Science and Technology Development Fund of Macau SAR (Grant No. 0002/2024/RDP).

\section*{Statements and Declarations}

\noindent\textbf{Competing Interests}

The authors have no competing interests to declare
that are relevant to the content of this article.

\noindent\textbf{Data Availability}

No datasets were generated or analyzed
during the current study.

\section*{Appendix}

\setcounter{equation}{0}
\renewcommand{\theequation}{A.\arabic{equation}}

\setcounter{definition}{0}
\renewcommand{\thedefinition}{A.\arabic{definition}}

\noindent\textbf{Proof of Proposition \ref{prop:Ytilde-bounded}:}
\begin{proof}
	Linearity follows from the linearity of the stochastic integral.
	The centeredness of \(Z_Y\) gives \(	\E[\widetilde Y(f)]
	=
	0. \)
	By \eqref{eq:Z-isometry} (see also \cite{Lindgren2012}),
	\begin{equation}\label{Y-BOUNDED}
		\begin{aligned}
			\E|\widetilde Y(f)|^2
			=
			\E\left|
			\int_{\mathbb R}
			\widehat f(\xi)\,Z_Y(\dd\xi)
			\right|^2
			=
			\int_{\mathbb R}
			|\widehat f(\xi)|^2\,\mu_Y(\dd\xi)
			=
			\int_{\mathbb R}
			|\widehat f(\xi)|^2S_Y(\xi)\,\dd\xi
			\leq
			\|S_Y\|_{L^\infty(\mathbb R)}
			\|\widehat f\|_{L^2(\mathbb R)}^2.
		\end{aligned}
	\end{equation}
	By Plancherel's theorem,
	$
	\|\widehat f\|_{L^2(\mathbb R)}
	=
	\|f\|_{L^2(\mathbb R)},
	$
	therefore
	\[
	\|\widetilde Y(f)\|_{L^2(\Omega)}^2
	\leq
	\|S_Y\|_{L^\infty(\mathbb R)}
	\|f\|_{L^2(\mathbb R)}^2.
	\]
\end{proof}

\noindent\textbf{Proof of Proposition \ref{prop:QY}:}
\begin{proof}
	By the Cauchy--Schwarz inequality and Plancherel's theorem,
	\[
	|Q_Y(f,g)|
	\leq
	\|S_Y\|_{L^\infty(\mathbb R)}
	\|\widehat f\|_{L^2(\mathbb R)}
	\|\widehat g\|_{L^2(\mathbb R)}		
	=
	\|S_Y\|_{L^\infty(\mathbb R)}
	\|f\|_{L^2(\mathbb R)}
	\|g\|_{L^2(\mathbb R)}.
	\]
	Moreover, since \(S_Y\geq0\) almost everywhere, \(	Q_Y(f,f)
	=
	\int_{\mathbb R}
	S_Y(\xi)
	|\widehat f(\xi)|^2\,\dd\xi
	\geq0.\)
\end{proof}
\noindent\textbf{Proof of Proposition \ref{prop:abstract-covariance-realization}:}
\begin{proof}
	By the Riesz representation theorem for bounded sesquilinear
	forms, there exists a unique bounded operator \(C_Q\)
	representing \(Q\). The positivity of \(Q\) implies that
	\(C_Q\) is positive and self-adjoint.
	Let \(W\) be Gaussian white noise on
	\(\overline{\operatorname{Ran}(C_Q^{1/2})}\),
	viewed as an isometric Gaussian random process
	(see \cite{FeichtingerHormann}).
	Define \(\mathcal X_Q(f)=W(C_Q^{1/2}f)\).
	Then \(\mathcal X_Q\) is centered and linear, with
	\[
	\E\bigl[\mathcal X_Q(f)\overline{\mathcal X_Q(g)}\bigr]
	=
	\inner{C_Q^{1/2}f}{C_Q^{1/2}g}_{L^2(\mathbb R)}
	=
	\inner{C_Qf}{g}_{L^2(\mathbb R)}
	=
	Q(f,g).
	\]
	Moreover,
	\(\|\mathcal X_Q(f)\|_{L^2(\Omega)}^2
	=Q(f,f)\leq M\|f\|_{L^2(\mathbb R)}^2\),
	so \(\mathcal X_Q\) is bounded.
\end{proof}
\noindent\textbf{Proof of Proposition \ref{second-order estimation}:}
\begin{proof}
	By \eqref{eq:covariance-operator-Ytilde}, \( \E|\widetilde Y(p_z)|^2
	=
	\inner{C_{\widetilde Y}p_z}{p_z}_{L^2(\R)}. \)
	Since $C_{\widetilde Y}$ is positive and bounded,
	\[
	\inner{C_{\widetilde Y}p_z}{p_z}_{L^2(\R)}
	\le
	\|C_{\widetilde Y}\|\|p_z\|_{L^2(\R)}^2
	=
	\frac{\pi\|C_{\widetilde Y}\|}{|\operatorname{Im}z|}.
	\]
\end{proof}

\noindent\textbf{Proof of Lemma \ref{lem:weighted-path-integrability}:}
\begin{proof}
	By second-order stationarity and the Cauchy--Schwarz
	inequality,
	\[
	\mathbb E|Y_t|
	\leq
	\bigl(\mathbb E|Y_t|^2\bigr)^{1/2}
	=
	\bigl(\mathbb E|Y_0|^2\bigr)^{1/2},
	\qquad t\in\mathbb R.
	\]
	Since $(t,\omega)\mapsto |Y_t(\omega)|/(1+t^2)$
	is nonnegative and jointly measurable, Tonelli's
	theorem gives
	\begin{equation}
		\mathbb E\left[
		\int_{\mathbb R}\frac{|Y_t|}{1+t^2}\,dt
		\right]
		=
		\int_{\mathbb R}\frac{\mathbb E|Y_t|}{1+t^2}\,dt
		\leq
		\bigl(\mathbb E|Y_0|^2\bigr)^{1/2}
		\int_{\mathbb R}\frac{dt}{1+t^2}
		=
		\pi\bigl(\mathbb E|Y_0|^2\bigr)^{1/2}
		<\infty.
	\end{equation}
	Thus the nonnegative random variable
	$\int_{\mathbb R}|Y_t|/(1+t^2)\,dt$
	is finite almost surely.
\end{proof}
\noindent\textbf{Proof of Proposition \ref{prop:cauchy-analyticity}:}
\begin{proof}
	We prove the assertion on the upper half-plane.
	By Lemma~\ref{lem:weighted-path-integrability}, the event
	\[
	A=
	\left\{
	\omega\in\Omega:
	\int_{\mathbb R}\frac{|Y_t(\omega)|}{1+t^2}\,dt<\infty
	\right\}
	\]
	has probability one.
	Recall that \({\widetilde Y}(f)=\int_{\mathbb R}\widehat f(\xi)\,Z(d\xi),
	~f\in L^2(\mathbb R),\)
	and $U_+(\omega,z)={\widetilde Y}(p_z)$.
	For $q\in L^1(\mathbb R)\cap L^2(\mathbb R)$,
	\eqref{eq:qY-path-integrable} in
	Remark~\ref{rem:sample-path-integrability} ensures the almost sure
	absolute convergence of the time-domain integral, while
	Remark~\ref{rem:spectral-bochner-interchange} justifies the
	interchange of $I$ with the Bochner integral. Thus,
	for almost every $\omega\in\Omega$,
	\begin{equation}\label{eq:Xq-time-domain}
		\begin{aligned}
			{\widetilde Y}(q)(\omega)
			&=\int_{\mathbb R}\widehat q(\xi)\,Z(d\xi)
			=I(\hat q)=\int_{\mathbb R}
			\left(
			\frac{1}{\sqrt{2\pi}}
			\int_{\mathbb R}q(t)e^{-it\xi}\,dt
			\right)Z(d\xi)\\
			&=\int_{\mathbb R}q(t)
			\left(
			\frac{1}{\sqrt{2\pi}}
			\int_{\mathbb R}e^{-it\xi}\,Z(d\xi)
			\right)dt
			=\int_{\mathbb R}q(t)Y_t(\omega)\,dt.
		\end{aligned}
	\end{equation}
	In particular,
	\begin{equation}\label{eq:Xq-pathwise}
		{\widetilde Y}(q)(\omega)
		=
		\int_{\mathbb R}q(t)Y_t(\omega)\,dt
		\qquad a.s.
	\end{equation}	
	To obtain a version defined simultaneously for all
	$z\in\mathbb C^+$, fix $a\in\mathbb C^+$.
	Since
	\[
	|p_z(t)-p_a(t)|
	\leq
	\frac{C_{z,a}}{1+t^2},
	\qquad t\in\mathbb R,
	\]
	we have
	$p_z-p_a\in L^1(\mathbb R)\cap L^2(\mathbb R)$.
	By \eqref{eq:Xq-pathwise},
	\[
	{\widetilde Y}(p_z-p_a)(\omega)
	=
	\int_{\mathbb R}
	Y_t(\omega)\bigl(p_z(t)-p_a(t)\bigr)\,dt
	\qquad\text{a.s.}
	\]
	We therefore choose the modification
	\begin{equation}\label{eq:Uplus-modification}
		U_+(\omega,z)
		=
		U_+(\omega,a)
		+
		\int_{\mathbb R}
		Y_t(\omega)\bigl(p_z(t)-p_a(t)\bigr)\,dt,
		\qquad \omega\in A,
	\end{equation}
	and set it equal to zero for $\omega\notin A$.
	The integral is absolutely convergent for every
	$\omega\in A$ and $z\in\mathbb C^+$.
	By \eqref{eq:Xq-pathwise}, this agrees almost surely
	with $\widetilde Y(p_z)$ for each fixed $z$.
	We now establish analyticity.
	Fix $\omega\in A$ and $z_0=x_0+iy_0\in\mathbb C^+$.
	For $0<|h|<y_0/2$, set
	\[
	q_h(t)
	=
	\frac{p_{z_0+h}(t)-p_{z_0}(t)}{h}
	=
	\frac{1}{(t-z_0-h)(t-z_0)}.
	\]
	By \eqref{eq:Uplus-modification},
	\begin{equation}
		\frac{U_+(\omega,z_0+h)-U_+(\omega,z_0)}{h}
		=\int_{\mathbb R}Y_t(\omega)q_h(t)\,dt
		=\int_{\mathbb R}
		\frac{Y_t(\omega)}
		{(t-z_0-h)(t-z_0)}\,dt.
	\end{equation}
	Since $|t-z_0|\geq y_0$ and $|h|<y_0/2$, \(	|t-z_0-h|
	\geq |t-z_0|-|h|
	\geq \frac12|t-z_0|. \)
	Moreover,
	\[
	\frac{1+t^2}{|t-z_0|^2}
	=
	\frac{1+t^2}{(t-x_0)^2+y_0^2}
	\longrightarrow 1
	\qquad\text{as }|t|\to\infty.
	\]
	Since this ratio is continuous on $\mathbb R$,
	it is bounded. Hence, for some $C_{z_0}>0$,
	\begin{equation}\label{eq:qh-domination}
		|Y_t(\omega)q_h(t)|
		\leq
		\frac{2|Y_t(\omega)|}{|t-z_0|^2}
		\leq
		2C_{z_0}\frac{|Y_t(\omega)|}{1+t^2}.
	\end{equation}
	The right-hand side is integrable for every
	$\omega\in A$ by
	Lemma~\ref{lem:weighted-path-integrability}
	and is independent of $h$.
	Since $q_h(t)\longrightarrow 1/(t-z_0)^2 \quad\text{as }h\to0$,
	the dominated convergence theorem gives
	\[
	\lim_{h\to0}
	\frac{U_+(\omega,z_0+h)-U_+(\omega,z_0)}{h}
	=
	\lim_{h\to0}
	\int_{\mathbb R}Y_t(\omega)q_h(t)\,dt
	=
	\int_{\mathbb R}
	\frac{Y_t(\omega)}{(t-z_0)^2}\,dt.
	\]
	The limiting integral is absolutely convergent, since
	\[
	\int_{\mathbb R}
	\frac{|Y_t(\omega)|}{|t-z_0|^2}\,dt
	\leq
	C_{z_0}
	\int_{\mathbb R}
	\frac{|Y_t(\omega)|}{1+t^2}\,dt
	<\infty.
	\]
	Since $z_0\in\mathbb C^+$ and $\omega\in A$
	were arbitrary, $U_+$ has analytic sample paths
	on $\mathbb C^+$ for every $\omega\in A$.
	The same argument applies to $U_-$ on
	$\mathbb C^-$ on the same event $A$.
	As $\mathbb P(A)=1$, both processes have
	analytic sample paths almost surely.
\end{proof}
\begin{remark}\label{rem:sample-path-integrability}
	
	Let $Y_t$ and $\widetilde Y$ be given by
	\eqref{eq:spectral-representation-Y} and \eqref{DEF-YTILDE},
	respectively. For every
	$q\in L^1(\mathbb R)\cap L^2(\mathbb R)$,
	the random variable $\widetilde Y(q)$ is well defined since
	$q\in L^2(\mathbb R)$.
	Moreover, by Tonelli's theorem, the Cauchy--Schwarz inequality,
	and second-order stationarity,
	\[
	\mathbb E
	\left[
	\int_{\mathbb R}|q(t)Y_t|\,dt
	\right]
	=
	\int_{\mathbb R}|q(t)|\,\mathbb E|Y_t|\,dt
	\leq
	\int_{\mathbb R}|q(t)|
	\bigl(\mathbb E|Y_t|^2\bigr)^{1/2}\,dt
	=
	\bigl(\mathbb E|Y_0|^2\bigr)^{1/2}
	\|q\|_{L^1(\mathbb R)}
	<\infty .
	\]
	Consequently,
	\begin{equation}\label{eq:qY-path-integrable}
		\int_{\mathbb R}|q(t)Y_t(\omega)|\,dt<\infty
		\qquad\text{for almost every }\omega\in\Omega.
	\end{equation}
	Hence \(	\int_{\mathbb R}q(t)Y_t(\omega)\,dt\)
	is an absolutely convergent sample-path integral for almost every
	$\omega\in\Omega$.

\end{remark}

\begin{remark}\label{rem:spectral-bochner-interchange}
	
	By Remark~\ref{rem:sample-path-integrability},
	the time-domain integral in \eqref{eq:Xq-time-domain}
	is well defined almost surely.
	We justify the interchange of the spectral and time integrals
	using the operator $I$ defined in \eqref{DEF-I}.
	For $q\in L^1(\mathbb R)\cap L^2(\mathbb R)$, set
	$F_q(t)(\xi)=q(t)e^{-it\xi}/\sqrt{2\pi}$.
	By the spectral representation of $Y_t$ and
	\eqref{eq:Z-isometry}, we have
	$\mathbb E|Y_t|^2=\mu_Y(\mathbb R)/(2\pi)$.
	Thus, $\mu_Y(\mathbb R)=2\pi\mathbb E|Y_0|^2<\infty$,
	and
	\[
	\|F_q(t)\|_{L^2(\mu_Y)}^2
	=
	\frac{|q(t)|^2}{2\pi}\mu_Y(\mathbb R)
	<
	\infty
	\]
	for almost every $t\in\mathbb R$.
	Since $\mu_Y(\mathbb R)<\infty$,
	the map $t\mapsto e^{-it\,\cdot}$ is continuous
	into $L^2(\mu_Y)$ by dominated convergence.
	Thus $F_q$ is strongly measurable, and
	\begin{equation}\label{eq:Fq-bochner-integrability}
		\int_{\mathbb R}
		\|F_q(t)\|_{L^2(\mu_Y)}\,dt
		=
		\frac{\mu_Y(\mathbb R)^{1/2}}{\sqrt{2\pi}}
		\int_{\mathbb R}|q(t)|\,dt
		=
		\bigl(\mathbb E|Y_0|^2\bigr)^{1/2}
		\|q\|_{L^1(\mathbb R)}
		<
		\infty.
	\end{equation}
	Therefore, $F_q$ is Bochner integrable in
	$L^2(\mu_Y)$, and by the Fourier transform,
	\begin{equation}\label{eq:Fq-fourier}
		\int_{\mathbb R}F_q(t)\,dt
		=
		\int_{\mathbb R}
		\frac{q(t)}{\sqrt{2\pi}}e^{-it\xi}\,dt
		=
		\widehat q.
	\end{equation}
	By \eqref{DEF-I} and the spectral representation
	of $Y_t$, we have $I(F_q(t))=q(t)Y_t$.
	Since $I$ is bounded and commutes with
	Bochner integrals \cite{vanNeerven2008},
	\eqref{Y-I-COMPOSITION} yields
	\begin{equation}\label{eq:spectral-bochner-interchange}
		\widetilde Y(q)
		=
		I(\widehat q)
		=
		I\left(
		\int_{\mathbb R}F_q(t)\,dt
		\right)
		=
		\int_{\mathbb R}I(F_q(t))\,dt
		=
		\int_{\mathbb R}q(t)Y_t\,dt
	\end{equation}
	in $L^2(\Omega)$.
	Consequently, \eqref{eq:Xq-pathwise} holds
	almost surely. Here $q\in L^2(\mathbb R)$
	ensures that $\widetilde Y(q)$ is defined,
	while $q\in L^1(\mathbb R)$ guarantees
	absolute integrability of its time-domain
	representation.

\end{remark}
\noindent\textbf{Proof of Lemma \ref{FT-cauchy}:}
\begin{proof}
	We use the Fourier transform convention
	$\widehat f(\xi)
	=
	\frac{1}{\sqrt{2\pi}}
	\int_{\mathbb R}f(t)e^{-it\xi}\dd t$.
	For $\xi\neq0$, consider
	$f(w)=\frac{e^{-iw\xi}}{w-z}$.
	By Jordan's lemma and the residue theorem,
	\[
	\int_{\mathbb R}
	\frac{e^{-it\xi}}{t-z}\dd t
	=
	\begin{cases}
		2\pi i e^{-iz\xi},
		& z\in\C^+,\ \xi<0,\\
		-2\pi i e^{-iz\xi},
		& z\in\C^-,\ \xi>0,\\
		0,
		& \text{otherwise}.
	\end{cases}
	\]
	The negative sign comes from the clockwise contour.
	Consequently,
	\[
	\widehat{p_z}(\xi)
	=
	\begin{cases}
		\sqrt{2\pi}\,i\,e^{-iz\xi}
		\1_{(-\infty,0)}(\xi),
		& z\in\C^+,\\[2pt]
		-\sqrt{2\pi}\,i\,e^{-iz\xi}
		\1_{(0,\infty)}(\xi),
		& z\in\C^-.
	\end{cases}
	\]
	The formulas hold almost everywhere, since
	the value at $\xi=0$ is irrelevant.
\end{proof}

\noindent\textbf{Proof of Theorem \ref{thm:two-sided-cauchy-complete}:}
\begin{proof}
	Suppose that \(f\in L^2(\mathbb R)\) is orthogonal to \(p_z\) for every
	\(z\in\mathbb C^+\cup\mathbb C^-\). By Plancherel's theorem and Lemma~\ref{FT-cauchy},
	for \(z=x+iy\in\mathbb C^+\), we have \(	0
	=
	\langle f,p_z\rangle_{L^2(\mathbb R)}
	=
	\langle \widehat f,\widehat{p_z}\rangle_{L^2(\mathbb R)}.\)
	Since \eqref{eq:FT-cauchy-plus}
	it follows that
	\[
	\int_{-\infty}^0
	\widehat f(\xi)e^{i\overline z\xi}\,d\xi
	=
	0,
	\qquad z\in\mathbb C^+.
	\]
	For fixed \(y>0\), this is the Fourier transform with respect to \(x\) of the \(L^1\)-function \(\widehat f(\xi)e^{y\xi}\mathbf 1_{(-\infty,0)}(\xi). \)
	Hence
	$
	\widehat f(\xi)=0
	~\text{for a.e. } \xi<0.
	$
	Similarly, using the formula
	\eqref{eq:FT-cauchy-minus}
	we obtain
	$
	\widehat f(\xi)=0
	~\text{for a.e. } \xi>0.
	$
	Therefore \(\widehat f=0\) almost everywhere on \(\mathbb R\), and hence
	\(f=0\). Thus the orthogonal complement of
	$
	\operatorname{span}\{p_z:z\in\mathbb C^+\cup\mathbb C^-\}
	$
	is trivial. Consequently,
	$
	\overline{\operatorname{span}}\{p_z:z\in\mathbb C^+\cup\mathbb C^-\}
	=
	L^2(\mathbb R).
	$
\end{proof}
\noindent\textbf{Proof of Lemma \ref{lem:weight-integrability}:}
\begin{proof}
	For $z=x+iy\in\mathbb C^+$ with $y>0$,
	\eqref{half-plane weight} gives
	\[
	I_\alpha
	=
	\int_0^\infty
	\int_{\mathbb R}
	y^{\alpha-1}
	\bigl((1+y)^2+x^2\bigr)^{-\alpha}
	\,dx\,dy.
	\]
	For fixed $y>0$, the substitution $x=(1+y)u$ gives
	\[
	\int_{\mathbb R}
	\bigl((1+y)^2+x^2\bigr)^{-\alpha}\,dx
	=
	(1+y)^{1-2\alpha}
	\int_{\mathbb R}(1+u^2)^{-\alpha}\,du.
	\]
	Set
	$c_\alpha=\int_{\mathbb R}(1+u^2)^{-\alpha}\,du<\infty$
	for $\alpha>1$. Hence,
	\[
	I_\alpha
	=
	c_\alpha
	\int_0^\infty
	y^{\alpha-1}(1+y)^{1-2\alpha}\,dy
	=
	c_\alpha B(\alpha,\alpha-1)
	<
	\infty.
	\]
	Indeed, the integrand behaves like $y^{\alpha-1}$
	as $y\to0$ and $y^{-\alpha}$ as $y\to\infty$.
	For the lower half-plane, the substitution $v=-y>0$
	reduces the corresponding integral to $I_\alpha$.
	Thus, both integrals are finite for $\alpha>1$.
\end{proof}

\noindent\textbf{Proof of Theorem \ref{THM-Bochner}:}
\begin{proof}
	By Proposition~\ref{prop:cauchy-analyticity},
	$U_+$ and $U_-$ have almost surely analytic
	sample paths on $\Cp$ and $\Cm$, respectively.
	Moreover, by Tonelli's theorem,
	Proposition~\ref{second-order estimation},
	and Lemma~\ref{lem:weight-integrability},
	\[
	\begin{aligned}
		\E\|U_+\|_{B_\alpha^2(\Cp)}^2
		&\leq
		\pi\|C_{\widetilde Y}\|
		\int_{\Cp}
		\frac{\rho_\alpha(z)}
		{\operatorname{Im}z}
		\,\dd A(z)
		=
		\pi\|C_{\widetilde Y}\|I_\alpha<\infty,
		\\
		\E\|U_-\|_{B_\alpha^2(\Cm)}^2
		&\leq
		\pi\|C_{\widetilde Y}\|
		\int_{\Cm}
		\frac{\rho_\alpha(z)}
		{|\operatorname{Im}z|}
		\,\dd A(z)
		=
		\pi\|C_{\widetilde Y}\|I_\alpha<\infty.
	\end{aligned}
	\]
	Together with the almost sure analyticity,
	these estimates imply that
	\[
	U_+\in L^2\bigl(\Omega,B_\alpha^2(\Cp)\bigr),
	\qquad
	U_-\in L^2\bigl(\Omega,B_\alpha^2(\Cm)\bigr).
	\]
	By the definition of the direct sum norm,
	\[
	\E\|U\|_{\cB}^2
	=
	\E\|U_+\|_{B_\alpha^2(\Cp)}^2
	+
	\E\|U_-\|_{B_\alpha^2(\Cm)}^2
	\leq
	2\pi\|C_{\widetilde Y}\|I_\alpha
	<
	\infty.
	\]
	Therefore,
	$
	U\in L^2(\Omega,\cB).
	$
\end{proof}

\noindent\textbf{Proof of Lemma \ref{lem:bergman-rkhs}:}
\begin{proof}
	Fix $\sigma\in\{+,-\}$ and $a\in\mathbb C^\sigma$.
	For any $0<r<|\operatorname{Im}a|$, the weight
	$\rho_\alpha$ is continuous and strictly positive
	on $\overline{D(a,r)}$. Hence, there exists
	$m_{a,r}>0$ such that
	$\rho_\alpha(z)\ge m_{a,r}$ for
	$z\in\overline{D(a,r)}$.	
	For $F\in B^2_\alpha(\mathbb C^\sigma)$,
	the function $|F_\sigma|^2$ is subharmonic.
	By the sub-mean value property,
	\[
	|F_\sigma(a)|^2
	\leq
	\frac{1}{\pi r^2}
	\int_{D(a,r)}|F_\sigma(z)|^2\,dA(z)\\
	\leq
	\frac{1}{\pi r^2m_{a,r}}
	\int_{D(a,r)}|F_\sigma(z)|^2\rho_\alpha(z)\,dA(z)
	\leq
	\frac{1}{\pi r^2m_{a,r}}
	\|F_\sigma\|_{B^2_\alpha(\mathbb C^\sigma)}^2.
	\]
	Setting
	$C_a:=\sup_{0<r<|\operatorname{Im}a|}
	\pi r^2m_{a,r}$,
	we have $0<C_a<\infty$ and
	$|F_\sigma(a)|^2\leq C_a^{-1}
	\|F_\sigma\|_{B^2_\alpha(\mathbb C^\sigma)}^2$.
	Thus the point evaluation functional
	$\Lambda_a:F_\sigma\mapsto F_\sigma(a)$ is bounded.
	By the Riesz representation theorem,
	there exists a unique
	$K_a^{\mathbb C^\sigma}
	\in B^2_\alpha(\mathbb C^\sigma)$ such that
	\[
	F(a)
	=
	\left\langle F,K_a^{\mathbb C^\sigma}\right\rangle_
	{B^2_\alpha(\mathbb C^\sigma)},
	\qquad
	F\in B^2_\alpha(\mathbb C^\sigma).
	\]
	Hence $B^2_\alpha(\mathbb C^\sigma)$ is a
	reproducing kernel Hilbert space.
\end{proof}

\noindent\textbf{Proof of Lemma \ref{lem:boundary-vanishing-half-plane-detailed}:}
\begin{proof}
	We prove the result for the upper half-plane;
	the lower half-plane case is analogous.
	For $z=x+iy\in\mathbb C^+$, the weight $\rho_\alpha(z)$
	is given by the upper half-plane case of \eqref{half-plane weight}.
	Consider the Cayley transform
	$\phi_+(w)=i\frac{1+w}{1-w}$,
	which maps $\mathbb D$ conformally onto $\mathbb C^+$.
	Since
	$\operatorname{Im}\phi_+(w)
	=\frac{1-|w|^2}{|1-w|^2}$,
	$\phi_+(w)+i=\frac{2i}{1-w}$,
	and
	$dA(\phi_+(w))
	=\frac{4}{|1-w|^4}\,dA(w)$,
	we obtain
	\begin{equation}\label{rho_phi_dA}
		\rho_\alpha(\phi_+(w))\,dA(\phi_+(w))
		=
		4^{1-\alpha}
		\frac{(1-|w|^2)^\alpha}{|1-w|^4}\,dA(w).
	\end{equation}
	Define
	\[
	(T_+F_+)(w)
	=
	2^{1-\alpha}
	\frac{F_+(\phi_+(w))}{(1-w)^2},
	\qquad
	F_+\in B^2_\alpha(\mathbb C^+),\quad w\in\mathbb D.
	\]
	Then $T_+F_+$ is holomorphic on $\mathbb D$.
	By the change of variables $z=\phi_+(w)$
	and \eqref{rho_phi_dA},
	\[
	\|F_+\|_{B^2_\alpha(\mathbb C^+)}^2
	=
	\int_{\mathbb C^+}
	|F_+(z)|^2\rho_\alpha(z)\,dA(z)
	=
	\int_{\mathbb D}
	4^{1-\alpha}
	\frac{|F_+(\phi_+(w))|^2}{|1-w|^4}
	(1-|w|^2)^\alpha\,dA(w)
	=
	\|T_+F_+\|_{A^2_\alpha(\mathbb D)}^2,
	\]
	where
	\(	A^2_\alpha(\mathbb D)
	=
	\left\{
	G\in\operatorname{Hol}(\mathbb D):
	\int_{\mathbb D}|G(w)|^2
	(1-|w|^2)^\alpha\,dA(w)<\infty
	\right\}.\)
	Moreover, $T_+$ is onto. Indeed, for any
	$G\in A^2_\alpha(\mathbb D)$, define
	\(	F_+(z)
	=
	2^{\alpha-1}
	\bigl(1-\phi_+^{-1}(z)\bigr)^2
	G(\phi_+^{-1}(z)),
	~ z\in\mathbb C^+.\)
	Then $T_+F_+=G$. Thus $T_+$ is a unitary
	isomorphism, and
	$\langle F_+,H_+\rangle_{B^2_\alpha(\mathbb C^+)}
	=
	\langle T_+F_+,T_+H_+\rangle_{A^2_\alpha(\mathbb D)}$.	
	By the boundary vanishing property of
	$A^2_\alpha(\mathbb D)$ established in
	\cite{QuDang2018}, for every
	$G\in A^2_\alpha(\mathbb D)$,
	\[
	\lim_{|b|\to1}
	\left|
	\left\langle G,e_b^{\mathbb D}
	\right\rangle_{A^2_\alpha(\mathbb D)}
	\right|
	=0,
	\qquad
	e_b^{\mathbb D}
	=
	\frac{K_b^{\mathbb D}}
	{\|K_b^{\mathbb D}\|_{A^2_\alpha(\mathbb D)}}.
	\]	
	We transfer this property to $\mathbb C^+$.
	For $a=x+iy\in\mathbb C^+$, let
	$b=\phi_+^{-1}(a)=\frac{a-i}{a+i}$.
	Then
	\[
	|b|^2
	=
	1-\frac{4y}{x^2+(y+1)^2},
	\]
	so $|b|\to1$ as $a$ approaches the real axis
	or tends to infinity in $\mathbb C^+$.
	
	For $F_+\in B^2_\alpha(\mathbb C^+)$,
	let $G=T_+F_+$ and
	$m_b:=2^{\alpha-1}(1-b)^2$.
	By the definition of $T_+$,
	$F_+(a)=m_bG(b)$.
	The reproducing properties and the unitarity
	of $T_+$ give
	\[
	\langle G,T_+K_a^{\mathbb C^+}
	\rangle_{A^2_\alpha(\mathbb D)}
	=F_+(a)=m_bG(b)
	=
	\langle G,\overline{m_b}K_b^{\mathbb D}
	\rangle_{A^2_\alpha(\mathbb D)}.
	\]
	Since this holds for every
	$G\in A^2_\alpha(\mathbb D)$, \(	T_+K_a^{\mathbb C^+}
	=
	\overline{m_b}K_b^{\mathbb D},
	~
	\|K_a^{\mathbb C^+}\|_{B^2_\alpha(\mathbb C^+)}
	=
	|m_b|\|K_b^{\mathbb D}\|_{A^2_\alpha(\mathbb D)}.\)
	Consequently,
	\(	T_+e_a^+
	=
	\frac{\overline{m_b}}{|m_b|}
	e_b^{\mathbb D}
	=
	\lambda_b e_b^{\mathbb D},
	\qquad |\lambda_b|=1.\)
	Therefore,
	\[
	\left|
	\left\langle F_+,e_a^+\right\rangle_
	{B^2_\alpha(\mathbb C^+)}
	\right|
	=
	\left|
	\left\langle T_+F_+,T_+e_a^+\right\rangle_
	{A^2_\alpha(\mathbb D)}
	\right|
	=
	\left|
	\left\langle G,e_b^{\mathbb D}\right\rangle_
	{A^2_\alpha(\mathbb D)}
	\right|
	\longrightarrow 0
	\]
	as $a\to\partial\mathbb C^+$, since $|b|\to1$.
	
	The lower half-plane case follows analogously
	from the Cayley transform
	$\phi_-(w)=-i\frac{1+w}{1-w}$,
	which induces the isometric isomorphism
	\(	(T_-F_-)(w)
	=
	2^{1-\alpha}
	\frac{F_-(\phi_-(w))}{(1-w)^2}\)
	from $B^2_\alpha(\mathbb C^-)$ onto
	$A^2_\alpha(\mathbb D)$.
	The same argument gives
	\(	\lim_{b\to\partial\mathbb C^-}
	\left|
	\left\langle F_-,e_b^-\right\rangle_
	{B^2_\alpha(\mathbb C^-)}
	\right|
	=0.\)
\end{proof}

\noindent\textbf{Proof of Corollary \ref{cor:boundary-vanishing}:}
\begin{proof}
	For \(a\in\mathbb C^+\), by \eqref{E+a}
	and the inner product on the Hilbert direct sum,
	\(
	\langle F,E_a^+\rangle_{\mathcal B_\alpha}
	=
	\langle F_+,e_a^+\rangle_{B_\alpha^2(\mathbb C^+)}.
	\)
	Hence, by Lemma~\ref{lem:boundary-vanishing-half-plane-detailed},
	\[
	\lim_{a\to\partial\mathbb C^+}
	\left|
	\left\langle F,E_a^+\right\rangle_{\mathcal B_\alpha}
	\right|
	=
	\lim_{a\to\partial\mathbb C^+}
	\left|
	\left\langle F_+,e_a^+\right\rangle_{B_\alpha^2(\mathbb C^+)}
	\right|
	=0.
	\]
	Similarly, for \(b\in\mathbb C^-\), by \eqref{E-b},
	\(
	\langle F,E_b^-\rangle_{\mathcal B_\alpha}
	=
	\langle F_-,e_b^-\rangle_{B_\alpha^2(\mathbb C^-)}.
	\)
	Again by Lemma~\ref{lem:boundary-vanishing-half-plane-detailed},
	\[
	\lim_{b\to\partial\mathbb C^-}
	\left|
	\left\langle F,E_b^-\right\rangle_{\mathcal B_\alpha}
	\right|
	=
	\lim_{b\to\partial\mathbb C^-}
	\left|
	\left\langle F_-,e_b^-\right\rangle_{B_\alpha^2(\mathbb C^-)}
	\right|
	=0.
	\]
\end{proof}
\noindent\textbf{Proof of Lemma \ref{lem:random-boundary-vanishing}:}
\begin{proof}
	Since \(U\in L^2(\Omega,\mathcal B_\alpha)\),
	we have \(U(\omega,\cdot)\in\mathcal B_\alpha\)
	for almost every \(\omega\in\Omega\).
	Write \(U\) as in \eqref{Cauchy transform}, with
	\(U_\pm(\omega,\cdot)
	=U(\omega,\cdot)|_{\mathbb C^\pm}
	\in B_\alpha^2(\mathbb C^\pm)\).
	For \(a\in\mathbb C^+\), by \eqref{E+a}
	and the inner product on \(\mathcal B_\alpha\),
	\[
	\left\langle U(\omega,\cdot),E_a^+
	\right\rangle_{\mathcal B_\alpha}
	=
	\left\langle U_+(\omega,\cdot),e_a^+
	\right\rangle_{B_\alpha^2(\mathbb C^+)}.
	\]
	Hence, by Corollary~\ref{cor:boundary-vanishing},
	for almost every \(\omega\in\Omega\),
	\[\left|
	\left\langle U(\omega,\cdot),E_a^+
	\right\rangle_{\mathcal B_\alpha}
	\right|^2
	\to 0
	\qquad\text{as }a\to\partial\mathbb C^+.\]
	Since \(\|E_a^+\|_{\mathcal B_\alpha}=1\),
	the Cauchy--Schwarz inequality gives \(\)
	\[
	\left|
	\left\langle U(\omega,\cdot),E_a^+
	\right\rangle_{\mathcal B_\alpha}
	\right|^2
	\leq
	\|U(\omega,\cdot)\|_{\mathcal B_\alpha}^2.
	\]
	The right-hand side is integrable because
	\(U\in L^2(\Omega,\mathcal B_\alpha)\).
	Thus, by the dominated convergence theorem,
	\(	\lim_{a\to\partial\mathbb C^+}
	\mathbb E
	\left|
	\left\langle U,E_a^+
	\right\rangle_{\mathcal B_\alpha}
	\right|^2
	=0.\)
	Similarly, for \(b\in\mathbb C^-\), by \eqref{E-b},
	\(	\langle U(\omega,\cdot),E_b^-\rangle_{\mathcal B_\alpha}
	=
	\langle U_-(\omega,\cdot),e_b^-\rangle_
	{B_\alpha^2(\mathbb C^-)}.\)
	Applying Corollary~\ref{cor:boundary-vanishing}
	and the same dominated convergence argument yields
	\(	\lim_{b\to\partial\mathbb C^-}
	\mathbb E
	\left|
	\left\langle U,E_b^-
	\right\rangle_{\mathcal B_\alpha}
	\right|^2
	=0.\)
\end{proof}

\noindent\textbf{Proof of Theorem \ref{thm:half-plane-orthogonality}:}
\begin{proof}
	For \(a\in\mathbb C^+\) and \(b\in\mathbb C^-\), by
	\eqref{eq:upper-lower-subspaces},
	\(	E_a^+\in B_\alpha^2(\mathbb C^+)\oplus\{0\},
	~
	E_b^-\in\{0\}\oplus B_\alpha^2(\mathbb C^-).\)
	Since \eqref{E+a} and \eqref{E-b},
	it follows from the definition of the inner product that
	\[
	\left\langle E_a^+,E_b^-\right\rangle_{\mathcal B_\alpha}
	=
	\left\langle e_a^+,0\right\rangle_{B_\alpha^2(\mathbb C^+)}
	+
	\left\langle 0,e_b^-\right\rangle_{B_\alpha^2(\mathbb C^-)}
	=0.
	\]
	Therefore, every element generated by the upper half-plane kernel family
	\(\{E_a^+:a\in\mathbb C^+\}\) is orthogonal to every element generated by
	the lower half-plane kernel family
	\(\{E_b^-:b\in\mathbb C^-\}\).
\end{proof}
\noindent\textbf{Proof of Theorem \ref{thm:SPOAFD-selection-existence}:}

\begin{proof}
	Continuity and \eqref{eq:SPOAFD-J-boundary} yield
	\eqref{eq:SPOAFD-maximizers}.
	Suppose \eqref{eq:SPOAFD-nonzero-residual} holds.
	If $J_n^+(a)=0$ for all
	$a\in\mathbb C^+$, then $G_{n,+}$ is
	orthogonal to all pre-orthogonalized
	candidates in the upper half-plane.
	Since it is also orthogonal to the previously
	selected elements, $G_{n,+}$ is orthogonal
	to the complete dictionary
	$\mathcal D_\alpha^+$, implying
	$G_{n,+}=0$ in
	$L^2(\Omega,B_\alpha^2(\mathbb C^+))$,
	a contradiction. Hence
	\eqref{eq:SPOAFD-positive-sup} holds.
	Similarly, if $G_{n,-}\neq0$ in
	$L^2(\Omega,B_\alpha^2(\mathbb C^-))$,
	the completeness of $\mathcal D_\alpha^-$
	implies
	\(	\sup_{b\in\mathbb C^-}J_n^-(b)>0.\)
	In both cases, the positive suprema are
	attained by \eqref{eq:SPOAFD-J-boundary}
	and continuity.
\end{proof}

\noindent\textbf{Proof of Theorem \ref{thm-spoafd-convergence}:}

\begin{proof}
	For every $N\geq1$, the systems
	$\{B_j^+\}_{j=1}^N$ and
	$\{B_\ell^-\}_{\ell=1}^N$
	are orthonormal in their respective
	weighted Bergman spaces.
	Thus, for almost every $\omega\in\Omega$,
	$U_{N,+}(\omega,\cdot)$ and
	$U_{N,-}(\omega,\cdot)$
	are the orthogonal projections of
	$U_+(\omega,\cdot)$ and
	$U_-(\omega,\cdot)$
	onto the corresponding finite-dimensional
	subspaces, respectively.
	By orthogonality and the zero-extension
	identification in Remark~\ref{zero-extended},
	we have
	\[
	\|U_+-U_{N,+}\|_{B_\alpha^2(\mathbb C^+)}^2
	=
	\|U_+\|_{B_\alpha^2(\mathbb C^+)}^2
	-
	\sum_{j=1}^N
	\left|
	\langle U,B_j^+\rangle_{\mathcal B_\alpha}
	\right|^2,
	\]
	and
	\[
	\|U_--U_{N,-}\|_{B_\alpha^2(\mathbb C^-)}^2
	=
	\|U_-\|_{B_\alpha^2(\mathbb C^-)}^2
	-
	\sum_{\ell=1}^N
	\left|
	\langle U,B_\ell^-\rangle_{\mathcal B_\alpha}
	\right|^2.
	\]
	Taking expectations leads to the two
	energy identities stated in the theorem.
	
	Since $\mathcal D_\alpha^+$ and
	$\mathcal D_\alpha^-$ are complete in their
	respective weighted Bergman spaces,
	together with lemma \ref{lem:boundary-vanishing-half-plane-detailed}
	and the  maximal selection rules,
	the assumptions of the SPOAFD convergence theorem
	\cite{Qian2022} are satisfied in both half-planes.
	Hence,
	\[
	\lim_{N\to\infty}
	\mathbb E
	\|U_+-U_{N,+}\|_{B_\alpha^2(\mathbb C^+)}^2
	=0
	\]
	and
	\[
	\lim_{N\to\infty}
	\mathbb E
	\|U_--U_{N,-}\|_{B_\alpha^2(\mathbb C^-)}^2
	=0.
	\]
	By the direct-sum structure of
	$\mathcal B_\alpha$ and
	\eqref{eq:G-U-U},
	\[
	\mathbb E\|U-U_N\|_{\mathcal B_\alpha}^2
	=
	\mathbb E
	\|U_+-U_{N,+}\|_{B_\alpha^2(\mathbb C^+)}^2
	+
	\mathbb E
	\|U_--U_{N,-}\|_{B_\alpha^2(\mathbb C^-)}^2.
	\]
	Therefore,
	\[
	\lim_{N\to\infty}
	\mathbb E\|U-U_N\|_{\mathcal B_\alpha}^2
	=0.
	\]
	This also establishes the expansions
	in \eqref{eq:U-expansion}
	in $L^2(\Omega,\mathcal B_\alpha)$.
\end{proof}

\noindent\textbf{Proof of Lemma \ref{PHI-F}:}
\begin{proof}
	For $z=x+iy\in\Cp$,
	$\|p_z\|_{L^2(\R)}=(\pi/y)^{1/2}$.
	By the Cauchy--Schwarz inequality,
	\[
	\begin{aligned}
		&\int_{\Cp}
		\|p_z\|_{L^2(\R)}
		|F_+(z)|\rho_\alpha(z)\,\dd A(z)
		=
		\sqrt{\pi}
		\int_{\Cp}
		|F_+(z)|y^{-1/2}\rho_\alpha(z)\,\dd A(z)\\
		&\quad\le
		\sqrt{\pi}
		\left(
		\int_{\Cp}
		|F_+(z)|^2\rho_\alpha(z)\,\dd A(z)
		\right)^{1/2}
		\left(
		\int_{\Cp}
		\frac{\rho_\alpha(z)}{y}\,\dd A(z)
		\right)^{1/2}
		=
		\sqrt{\pi I_\alpha}\,
		\|F_+\|_{B_\alpha^2(\Cp)}<\infty.
	\end{aligned}
	\]
	The same estimate holds on $\Cm$ for $F_-$.
	By Minkowski's integral inequality, the integrals
	in \eqref{Phi-F} converge absolutely for almost every
	$t\in\R$. Moreover,
	\[
	\|\Phi_F^+\|_{L^2(\R)}
	\le
	\int_{\Cp}
	\|p_z\overline{F_+(z)}\rho_\alpha(z)\|_{L^2(\R)}
	\,\dd A(z)
	=
	\int_{\Cp}
	\|p_z\|_{L^2(\R)}
	|F_+(z)|\rho_\alpha(z)\,\dd A(z)
	\le
	\sqrt{\pi I_\alpha}\,
	\|F_+\|_{B_\alpha^2(\Cp)}.
	\]
	Similarly,
	$\|\Phi_F^-\|_{L^2(\R)}
	\le\sqrt{\pi I_\alpha}\,
	\|F_-\|_{B_\alpha^2(\Cm)}$.
	Consequently,
	\[
	\|\Phi_F\|_{L^2(\R)}
	\le
	\|\Phi_F^+\|_{L^2(\R)}
	+
	\|\Phi_F^-\|_{L^2(\R)}
	\le
	\sqrt{\pi I_\alpha}
	\left(
	\|F_+\|_{B_\alpha^2(\Cp)}
	+
	\|F_-\|_{B_\alpha^2(\Cm)}
	\right)
	\le
	\sqrt{2\pi I_\alpha}\,\|F\|_{\cB}.
	\]
	
	Furthermore, the maps
	$z\mapsto p_z\overline{F_\pm(z)}\rho_\alpha(z)$
	are continuous into $L^2(\R)$ and have integrable
	norms. Hence their Bochner integrals are well defined
	and agree with $\Phi_F^\pm$ in $L^2(\R)$
	\cite{vanNeerven2008}.
	Since $\widetilde Y:L^2(\R)\to L^2(\Omega)$
	is bounded and linear, it commutes with
	Bochner integration \cite{vanNeerven2008}.
	Thus, in $L^2(\Omega)$,
	\[
	\begin{aligned}
		\widetilde Y(\Phi_F)
		&=
		\int_{\Cp}
		\widetilde Y(p_z)\overline{F_+(z)}
		\rho_\alpha(z)\,\dd A(z)
		+
		\int_{\Cm}
		\widetilde Y(p_z)\overline{F_-(z)}
		\rho_\alpha(z)\,\dd A(z)\\
		&=
		\int_{\Cp}
		U_+(z)\overline{F_+(z)}
		\rho_\alpha(z)\,\dd A(z)
		+
		\int_{\Cm}
		U_-(z)\overline{F_-(z)}
		\rho_\alpha(z)\,\dd A(z)\\
		&=
		\inner{U}{F}_{\cB}.
	\end{aligned}
	\]
\end{proof}
\noindent\textbf{Proof of Theorem \ref{thm:PHI-F}:}
\begin{proof}
	By \eqref{inner-U-B} and \eqref{eq:covariance-operator-Ytilde}, we obtain
	\(	\E\left|\inner{U}{F}_{\cB}\right|^2
	=
	\E|{\widetilde Y}(\Phi_F)|^2
	=
	\inner{C_{\widetilde Y}\Phi_F}{\Phi_F}_{L^2(\R)}.\)
	Since \(C_{\widetilde Y}\) is the Fourier multiplier with symbol \(S_Y\), Plancherel's theorem gives
	\[
	\inner{C_{\widetilde Y}\Phi_F}{\Phi_F}_{L^2(\R)}
	=
	\int_\R
	S_Y(\xi)|\widehat{\Phi_F}(\xi)|^2\,\dd\xi.
	\]
\end{proof}

\noindent\textbf{Proof of Theorem \ref{thm:SPOAFD-weak-recovery}:}

\begin{proof}
	For $t\in\mathbb R$ and $y>0$, let
	$a=t+iy\in\mathbb C^+$ and
	$\bar a=t-iy\in\mathbb C^-$, let  $V_y$ and $V_{N,y}$  be defined as in
	\eqref{eq:Poisson-regularization}
	and\eqref{eq:SPOAFD-Poisson-approximation} , respectively.
	We first show that, for each fixed $y>0$,
	\begin{equation}\label{eq:SPOAFD-test-convergence}
		\lim_{N\to\infty}
		\mathbb E\left|
		\int_{\mathbb R}
		\gamma(t)\bigl(V_{N,y}(t)-V_y(t)\bigr)\,dt
		\right|^2=0.
	\end{equation}
	Since $\gamma\in C_c^\infty(\mathbb R)$, there exists
	$R>0$ such that
	$\operatorname{supp}\gamma\subset[-R,R]$.
	Let $K_{a}^{\mathbb C^+}$ and $K_{\bar a}^{\mathbb C^-}$
	denote the reproducing kernels of
	$B_\alpha^2(\mathbb C^+)$ and
	$B_\alpha^2(\mathbb C^-)$, respectively.
	For $F\in\mathcal B_\alpha$, the reproducing
	property and the Cauchy--Schwarz inequality give
	\[
	|F_+(z)|
	=
	\left|
	\left\langle
	F_+,K_a^{\mathbb C^+}
	\right\rangle_{B_\alpha^2(\mathbb C^+)}
	\right|
	\leq
	\|F_+\|_{B_\alpha^2(\mathbb C^+)}
	\|K_a^{\mathbb C^+}\|_{B_\alpha^2(\mathbb C^+)}.
	\]
	Similarly, \(	|F_-(z)|
	\leq
	\|F_-\|_{B_\alpha^2(\mathbb C^-)}
	\|K_{\bar a}^{\mathbb C^-}\|_{B_\alpha^2(\mathbb C^-)},
	~ z\in\mathbb C^-.\)
	For fixed $y>0$, define
	\[
	M_{R}
	=
	\sup_{|t|\leq R}
	\Big(
	\|K_{a}^{\mathbb C^+}\|_{B_\alpha^2(\mathbb C^+)}
	+
	\|K_{\bar a}^{\mathbb C^-}\|_{B_\alpha^2(\mathbb C^-)}
	\Big).
	\]
	Since the norms of the reproducing kernels are
	locally bounded in their respective half-planes,
	we have $M_{R}<\infty$.
	Applying the preceding estimates to $U_N-U$,
	we obtain
	\[
	\begin{aligned}
		&\left|
		\int_{\mathbb R}
		\gamma(t)\bigl(V_{N,y}(t)-V_y(t)\bigr)\,dt
		\right|
		\leq
		\frac{1}{2\pi}
		\int_{-R}^{R}|\gamma(t)|
		\Big(
		|U_{N,+}(\omega,a)-U_+(\omega,a)|
		+
		|U_{N,-}(\omega,\bar a)-U_-(\omega,\bar a)|
		\Big)\,dt\\
		&\leq
		\frac{\|U_N-U\|_{\mathcal B_\alpha}}{2\pi}
		\int_{-R}^{R}|\gamma(t)|
		\Big(
		\|K_{a}^{\mathbb C^+}\|_{B_\alpha^2(\mathbb C^+)}
		+
		\|K_{\bar a}^{\mathbb C^-}\|_{B_\alpha^2(\mathbb C^-)}
		\Big)\,dt
		\leq
		\frac{M_{R}}{2\pi}
		\|\gamma\|_{L^1(\mathbb R)}
		\|U_N-U\|_{\mathcal B_\alpha}.
	\end{aligned}
	\]
	The preceding estimate also ensures that
	the test integral is absolutely convergent
	for almost every $\omega$.
	Taking squares and expectations gives
	\[
	\mathbb E\left|
	\int_{\mathbb R}
	\gamma(t)\bigl(V_{N,y}(t)-V_y(t)\bigr)\,dt
	\right|^2
	\leq
	\frac{M_{R}^2}{4\pi^2}
	\|\gamma\|_{L^1(\mathbb R)}^2
	\mathbb E\|U_N-U\|_{\mathcal B_\alpha}^2.
	\]
	By the mean-square convergence of the SPOAFD
	expansion, the right-hand side tends to zero
	as $N\to\infty$. Hence
	\eqref{eq:SPOAFD-test-convergence} follows.
	
	Next, recall that
	\[
	U_\pm(\omega,z)=\widetilde Y(p_z),
	\qquad
	p_z(s)=\frac{1}{s-z},
	\qquad z\in\mathbb C^\pm.
	\]
	For $t,s\in\mathbb R$, $y>0$, $ a=t+iy\in\mathbb C^+$ and
	$\bar a=t-iy\in\mathbb C^-$
	\[
	\frac{
		p_{a}(s)-p_{\bar a}(s)
	}{2\pi i}
	=
	\frac{1}{2\pi i}
	\left(
	\frac{1}{s-t-iy}
	-
	\frac{1}{s-t+iy}
	\right)
	=
	\frac{1}{\pi}
	\frac{y}{(s-t)^2+y^2}
	=
	P_y(t-s).
	\]
	By the linearity of $\widetilde Y$, it follows that
	\[
	V_y(\omega,t)
	=
	\frac{
		\widetilde Y(p_{a})
		-
		\widetilde Y(p_{\bar a})
	}{2\pi i}
	=
	\widetilde Y(P_y(t-\cdot))(\omega)
	\]
	in $L^2(\Omega)$.
	For each fixed $y>0$, the almost sure analyticity
	of $U_\pm$ implies that
	$t\mapsto V_y(\omega,t)$ is continuous
	for almost every $\omega$.
	Since $\gamma\in C_c^\infty(\mathbb R)$ with
	$\operatorname{supp}\gamma\subset[-R,R]$, we have
	\[
	\int_{\mathbb R}|\gamma(t)V_y(\omega,t)|\,dt
	\leq
	\|\gamma\|_{L^1(\mathbb R)}
	\sup_{|t|\leq R}|V_y(\omega,t)|
	<\infty
	\]
	for almost every $\omega$.
	Moreover, the map
	$t\mapsto \gamma(t)P_y(t-\cdot)$
	is Bochner integrable in $L^2(\mathbb R)$, since
	\[
	\int_{\mathbb R}
	|\gamma(t)|\|P_y(t-\cdot)\|_{L^2(\mathbb R)}\,dt
	=
	\|P_y\|_{L^2(\mathbb R)}
	\|\gamma\|_{L^1(\mathbb R)}
	<\infty.
	\]
	By the boundedness of $\widetilde Y$ and
	the Bochner integral interchange theorem \cite{vanNeerven2008},
	together with
	$V_y(\omega,t)=\widetilde Y(P_y(t-\cdot))(\omega)$
	in $L^2(\Omega)$, we obtain
	\[
	\int_{\mathbb R}\gamma(t)V_y(\omega,t)\,dt
	=
	\widetilde Y\left(
	\int_{\mathbb R}
	\gamma(t)P_y(t-\cdot)\,dt
	\right)(\omega)
	=
	\widetilde Y(P_y*\gamma)(\omega),
	\]
	for almost every $\omega$, where the last
	equality follows from the symmetry of
	the Poisson kernel.
	
	On the other hand, since
	$\gamma\in C_c^\infty(\mathbb R)
	\subset L^1(\mathbb R)\cap L^2(\mathbb R)$,
	Remark~\ref{rem:spectral-bochner-interchange}
	ensures that for any $\omega$, the time-domain integral
	$\int_{\mathbb R}\gamma(t)Y_t(\omega)\,dt$
	is well defined almost surely.
	Moreover, by
	\eqref{eq:spectral-bochner-interchange},
	\[
	\int_{\mathbb R}\gamma(t)Y_t(\omega)\,dt
	=
	\int_{\mathbb R}
	\left(
	\frac{1}{\sqrt{2\pi}}
	\int_{\mathbb R}
	\gamma(t)e^{-it\xi}\,dt
	\right)Z_Y(d\xi)
	=
	\int_{\mathbb R}
	\widehat\gamma(\xi)\,Z_Y(d\xi)
	=
	\widetilde Y(\gamma)(\omega)
	\]
	in $L^2(\Omega)$.
	Combining the preceding identities, we obtain
	\begin{equation}\label{eq:SPOAFD-Poisson-test-identity}
		\int_{\mathbb R}
		\gamma(t)\bigl(V_y(\omega,t)-Y_t(\omega)\bigr)\,dt
		=
		\widetilde Y(P_y*\gamma)(\omega)-\widetilde Y(\gamma)(\omega)
	\end{equation}
	in $L^2(\Omega)$.
	By the linearity of $\widetilde Y$ and the
	boundedness estimate \eqref{eq:Ytilde-bound},
	\[
	\mathbb E
	\left|
	\widetilde Y(P_y*\gamma)(\omega)-\widetilde Y(\gamma)(\omega)
	\right|^2
	=
	\mathbb E
	\left|
	\widetilde Y(P_y*\gamma-\gamma)(\omega)
	\right|^2
	\leq
	\|S_Y\|_{L^\infty(\mathbb R)}
	\|P_y*\gamma-\gamma\|_{L^2(\mathbb R)}^2.
	\]
	Since the Poisson kernels form an approximate identity
	in $L^2(\mathbb R)$, \(	\lim_{y\to0^+}
	\|P_y*\gamma-\gamma\|_{L^2(\mathbb R)}=0.\)
	It follows that
	\begin{equation}\label{eq:SPOAFD-Poisson-weak-limit}
		\lim_{y\to0^+}
		\mathbb E\left|
		\int_{\mathbb R}
		\gamma(t)\bigl(V_y(\omega,t)-Y_t\bigr(\omega))\,dt
		\right|^2=0.
	\end{equation}
	For every fixed $y>0$, the reverse triangle
	inequality in $L^2(\Omega)$ gives
	\[
	\begin{aligned}
		&\left|
		\left(
		\mathbb E\left|
		\int_{\mathbb R}
		\gamma(t)(V_{N,y}(\omega,t)-Y_t(\omega))\,dt
		\right|^2
		\right)^{1/2}
		-
		\left(
		\mathbb E\left|
		\int_{\mathbb R}
		\gamma(t)(V_y(\omega,t)-Y_t(\omega))\,dt
		\right|^2
		\right)^{1/2}
		\right|\\
		&\qquad\leq
		\left(
		\mathbb E\left|
		\int_{\mathbb R}
		\gamma(t)(V_{N,y}(\omega,t)-V_y(\omega,t))\,dt
		\right|^2
		\right)^{1/2}
		\longrightarrow 0
	\end{aligned}
	\]
	as $N\to\infty$, by
	\eqref{eq:SPOAFD-test-convergence}.
	Consequently,
	\[
	\lim_{N\to\infty}
	\mathbb E\left|
	\int_{\mathbb R}
	\gamma(t)(V_{N,y}(\omega,t)-Y_t(\omega))\,dt
	\right|^2
	=
	\mathbb E\left|
	\int_{\mathbb R}
	\gamma(t)(V_y(\omega,t)-Y_t(\omega))\,dt
	\right|^2.
	\]
	Taking $y\to0^+$ and applying
	\eqref{eq:SPOAFD-Poisson-weak-limit}, we obtain
	\[
	\lim_{y\to0^+}\lim_{N\to\infty}
	\mathbb E\left|
	\int_{\mathbb R}
	\gamma(t)\bigl(V_{N,y}(\omega,t)-Y_t(\omega)\bigr)\,dt
	\right|^2=0,
	\]
	which is precisely
	\eqref{eq:SPOAFD-weak-recovery}.
\end{proof}

	

\end{document}